\documentclass{conm-p-l}

\usepackage{}

\newtheorem{theorem}{Theorem}[section]

\newtheorem{proposition}[theorem]{Proposition}

\theoremstyle{definition}

\theoremstyle{remark}

\numberwithin{equation}{section}

\newcommand{\msn}{\medskip\noindent}

\newcommand{\spann}{\operatorname{span}}
\newcommand{\supp}{\operatorname{supp}}
\newcommand{\R}{\mathbb{R}}
\newcommand{\Rn}{\R^n}
\newcommand{\Rnm}{\R^{n\times m}}
\newcommand{\Rp}{\R_+}
\newcommand{\Rpn}{\Rp^n}
\newcommand{\Rpm}{\Rp^m}
\newcommand{\Rpnn}{\Rp^{n\times n}}
\newcommand{\Rpnm}{\Rp^{n\times m}}
\newcommand{\Rpnml}{\Rp^{n\times m_1}}
\newcommand{\Rpnmk}{\Rp^{n\times m_k}}

\newcommand{\kroneker}{\delta}
\newcommand{\perm}{\operatorname{per}}
\newcommand{\per}{\operatorname{per}}
\newcommand{\treug}{\bigtriangledown}
\newcommand{\adjug}{\operatorname{adj}}
\newcommand{\hilb}{\operatorname{H}}
\newcommand{\Hilb}{\operatorname{H}}

\newcommand{\halfs}{{\mathcal H}}
\newcommand{\sle}{\operatorname{sl}}

\newcommand{\Lin}{\operatorname{L}}
\newcommand{\cS}{{\mathcal S}}

\begin{document}

\title[Multiorder, Kleene stars and Cyclic Projectors]{Multiorder, Kleene Stars and Cyclic Projectors in the Geometry of Max
Cones}


%    author one information
\author[S. Sergeev]{Serge\u{\i} Sergeev}
\address{Serge\u{\i} Sergeev, University of Birmingham, School of Mathematics, B15 2TT Edgbaston, Birmingham, UK}
%\curraddr{}
\email{sergeevs@maths.bham.ac.uk}
\thanks{This work is supported by the EPSRC grant RRAH12809, the RFBR grant 08-01-00601 and the joint
RFBR/CNRS grant 05-01-02807}

%%    author two information
%\author{}
%\address{}
%\curraddr{}
%\email{}
%\thanks{}

\date{}

\begin{abstract}
Max cones are the subsets of the nonnegative orthant $\R_+^n$ of the
$n$-dimensional real space $\R^n$ closed under scalar multiplication
and componentwise maximisation. Their study is motivated by some practical applications
which arise in discrete event systems, optimal scheduling and modelling of
synchronization problems in multiprocessor interactive systems.
We investigate the geometry of max cones,
concerning the role of the multiorder principle, the Kleene stars, and the cyclic
projectors.

The multiorder principle is closely related to the set covering conditions in max algebra,
and gives rise to important analogues of some theorems
of convex geometry. We show that, in particular,
this principle leads to a convenient representation of
certain nonlinear projectors onto max cones.
% and that it is important for the geometry of
%the cellular decomposition and cyclic projectors.

The Kleene stars are fundamental in max algebra since they accumulate weights of
optimal paths and yield generators
for max-algebraic eigenspaces of matrices. We examine the role of their column spans
called Kleene cones, as building blocks in the Develin-Sturmfels cellular decomposition.
Further we show that the cellular decomposition gives rise to new max-algebraic
objects which we call row and column Kleene stars. We relate these objects to the max-algebraic
pseudoinverses of matrices and to tropical versions of the colourful Carath\'{e}odory theorem.

The cyclic projectors are
specific nonlinear operators which lead to the so-called
alternating method for finding a solution to homogeneous two-sided systems of
max-linear equations. We generalize the alternating method
to the case of homogeneous multi-sided systems,
and we give a proof, which uses the cellular decomposition idea, that the
alternating method converges in a finite number of iterations
to a positive solution of a multi-sided
system if a positive solution exists.
We also present new bounds on the number of iterations of the
alternating method, expressed in terms of the Hilbert projective distance between max cones.
\end{abstract}

\maketitle

\newpage
\section{Introduction}

The nonnegative orthant $\Rpn$ of the $n$-dimensional real space $\Rn$ can be viewed
as an $n$-dimensional free semimodule over the
max-times semiring, which is the set of nonnegative numbers $\Rp$ equipped
with the operations of 'addition' $a\oplus b:=\max(a,b)$ and the
ordinary multiplication $a\otimes b:=a\times b$. The max-times semiring is denoted by
$\R_{\max,\times}=(\Rp,\oplus=\max,\otimes=\times)$.
Zero and unity of the semiring
coincide with the usual $0$ and $1$. For instance, in this semiring
$2\otimes 3=6$ and $2\oplus 3=3$. Subsemimodules of $\Rpn=\R_{\max,\times}^n$ are the
subsets of $\Rpn$ closed
under the componentwise maximization $\oplus$, and the usual multiplication by nonnegative
scalars. These subsemimodules will be called {\em max cones}, due to their obvious analogy
with convex cones. In a very important special case, max cones can
indeed be convex cones, but in general they are not convex, i.e., not stable under the usual componentwise
addition.

By max algebra we understand linear algebra over the semiring
$\R_{\max,\times}$, extending the $\max,\times$ arithmetic to
nonnegative matrices and vectors in the usual way. For instance, if
$A=(a_{ij})$ and $B=(b_{ij})$ are two matrices of appropriate sizes,
then $(A\oplus B)_{ij}=a_{ij}\oplus b_{ij}$, or $(A\otimes
B)_{ij}=\bigoplus_k a_{ik} b_{kj}$. The iterated product $A\otimes
A\otimes...\otimes A$ in which the symbol
$A$ appears $k$ times will be denoted by $%
A^{k}$. We assume that $A^0:=I$, the unit matrix. The sets like $\{1,\ldots,m\}$
or $\{1,\ldots,n\}$ will be denoted by $[m]$ or $[n]$ respectively,
and for a set of indices $M$, the number of elements in $M$
will be denoted by $|M|$.

The idempotency of addition $a\oplus a=a$ and the lack of
subtraction are important features of max algebra that make it
different from the nonnegative linear algebra.

Max algebra has been known for some time, and we mention here the
pioneering works of Cuninghame-Green \cite{CG-60,CG:79}, Yoeli
\cite{Yoe-61}, Vorobyev \cite{Vor-67}, Carr\'{e} \cite{Car-71}, Gondran and Minoux \cite{GM-84},
K. Zimmermann \cite{Zim-76}, and U. Zimmermann \cite{UZim:81},
among many others. Max algebra is often presented in the
settings which seem to be different from $\R_{\max,\times}$, namely,
over semirings $\R_{\max,+}=(\R\cup\{-\infty\},
\oplus=\max,\otimes=+)$ (max-plus semiring),
$\R_{\min,+}=(\R\cup\{+\infty\}, \oplus=\min,\otimes=+)$ (tropical
or min-plus semiring), or most exotically
$\R_{\min,\times}=(\R_+\cup\{+\infty\},\oplus=\min,\otimes=\times)$
(min-times semiring). All these semirings are isomorphic to each
other and to $\R_{\max,\times}$. Max algebra has important practical applications which
arise in discrete event systems and scheduling problems
\cite{BCOQ,CG:79,GauThesis}, and in modelling of synchronization problems in
multiprocessor interactive systems \cite{BA-08}. 

More generally, max algebra can be seen as a branch of tropical
mathematics, which is a rapidly developing field with applications
in mathematical physics, optimal control, algebraic geometry and
other research areas. See \cite{Lit-07<Intro>} for a recent survey,
and \cite{LM:05,LMS:07} for recent collections of papers.

The similarity between max cones and convex cones was understood in
the very beginning by Vorobyev \cite{Vor-67}, who used the name
'extremally convex cones' (instead of semimodules or spaces).
K. Zimmermann \cite{Zim-77} defined extremally convex sets, or
tropically/max-plus convex sets as it would be called now, and
proved a separation theorem of a point from a closed convex set.
This theorem was generalized and more transparent proofs were given
by Samborski\u{\i} and Shpiz \cite{SS-92}, Litvinov et al.
\cite{LMS-01}, Cohen et al. \cite{CGQ-04,CGQS-05}, and also Develin
and Sturmfels \cite{DS-04}, Joswig \cite{Jos-05}. We note that the
separation theorem of a point from a closed max cone, given below as
Theorem \ref{separ0}, is essentially the same result. In the
ordinary convex geometry, separation of a point from a convex set
easily leads to the separation of two convex sets from each other.
However, analogous statements for max cones arise differently
and are related to the investigation of certain nonlinear projectors onto max cones,
and their compositions called cyclic projectors, see Gaubert and Sergeev \cite{GS-07} and
Theorems \ref{specrad=dh} and \ref{separ} below. Remarkably, these
cyclic projectors also appear in the study of two-sided max-linear systems
of equations, see Cuninghame-Green and Butkovi\v{c} \cite{CGB}, and
lead to a pseudopolynomial method for finding solutions to such systems.
This will be discussed in the last section of the paper. We also note here that cyclic projectors
are special case of the multiplicative version of the 
min-max functions studied in \cite{CZ-02, CTGG-99, OP-97}.

The geometry of max cones can be thought of as a special case of the
multiorder convexity, a concept introduced by Mart\'{\i}nez-Legaz
and Singer \cite{MLS-91}. Although this idea was made explicit only
recently in a work by Ni\c{t}ic\u{a} and Singer \cite{NS-07b},
it is closely related to the set-covering conditions for $A\otimes
x=b$ systems in max algebra \cite{BCOQ,CG:79,Vor-67}. The multiorder
principle, see Propositions \ref{p:multiorder} and \ref{ex-minima}
below, leads to easy proofs of many statements concerning
generators, extremals and bases of max cones, see Butkovi\v{c} et
al.\cite{BSS-07}, including the tropical Carath\'{e}odory theorem,
and Minkowski's theorem about extremals of closed cones (also
Gaubert and Katz \cite{GK-07}). The multiorder principle is also
important for the tropical convexity approach, meaning works of
Develin, Sturmfels, Joswig, Yu et al. \cite{BY-06,DS-04,Jos-05},
since it describes max cones as intersections of staircases, and
their extremals as elements of bases of monomial ideals.
%See also the monograph by Miller and Sturmfels \cite{MS:06}.

Yet another approach to the geometry of max cones, though strongly related to
the previous one, is to represent max cones as cellular complexes, or, roughly
speaking, as unions of ordinary convex cones. This approach was put forward by
Develin and Sturmfels \cite{DS-04}, and called cellular decomposition.
The atoms of this decomposition are well-known to specialists in convex geometry and
combinatorics, see
Joswig and Kulas \cite{JK-08} for more details.
As it was noticed in \cite{Ser-07}, these atoms are column spans of uniquely
defined Kleene stars, a fundamental concept in max algebra.

The aim of the present paper is to bring together some geometric and
algebraic ideas discussed above. Section 2 discusses the multiorder principle
and related results.
In particular, we show that this principle leads to a convenient new representation
of the nonlinear projectors mentioned above.
In Section 3 we recall the concept of Kleene stars and examine the role of their column spans
called Kleene cones as building blocks in the Develin-Sturmfels cellular decomposition.
Further we show that, in turn, the cellular decomposition gives rise to new max-algebraic objects which we call
row and column Kleene stars. We relate these new concepts to the max-algebraic pseudoinverses
of matrices and to tropical versions of the colourful Carath\'{e}odory theorem.
In Section 4 we generalize the alternating method of
Cuninghame-Green and Butkovi\v{c} \cite{CGB} to the case of
multisided systems $A^{(1)}\otimes x^1=\ldots=A^{(k)}\otimes x^k$.
We give a proof, based on the cellular decomposition idea,
that if the system has a positive
solution, then the method converges to a positive solution in a
finite number of steps. We also present new bounds for the number of
iterations in the max-plus integer case, and in the general case
when there are no solutions, in terms of the Hilbert projective distance between max cones.

\section{The role of multiorder}

\subsection{Generators, bases and extremals of max cones}
Let $S\subseteq\Rpn$. A vector $u\in\Rpn$ is called a {\em max combination} of $S$
if
\begin{equation}
\label{def:maxcomb}
u=\bigoplus_{v\in S} \lambda_v v,\ \lambda_v\in\Rp,
\end{equation}
where only a finite number of $\lambda_v$ are nonzero. The set of all
max combinations \eqref{def:maxcomb} of $S$ will be denoted by
$\spann(S)$. Evidently, $\spann(S)$ is a max cone. If $\spann(S)=V$,
then we call $S$ a set of {\em generators} for $V$ and say that $V$
is {\em generated}, or {\em spanned}, by $S$. In particular, the set
of all max combinations of columns of a matrix $A$ will be denoted
by $\spann(A)$ and called the {\em column span} of $A$. If none of
the elements of a generating set $S$ of a max cone $V$ can be
expressed as a max combination of other elements, then $S$ is called
a {\em (weak) basis} of $V$.

A vector $v\in V$ is called an {\em extremal} of $V$, if
\begin{equation*}
\label{def:extr}
v=u\oplus w,\, u,w\in V\Rightarrow \text{$v=u$ or $v=w$}.
\end{equation*}
Extremals are analogous to extremal rays of convex cones. If $v$ is an extremal of $V$ and
$\lambda>0$, then $\lambda v$ is also an extremal.

For all $i=1,\ldots,n$ define the following preorder relation.

\begin{equation*}
\label{def:preorder}
u\leq_j v\Leftrightarrow uu_j^{-1}\leq vv_j^{-1},\ u_j\neq 0,\ v_j\neq 0.
\end{equation*}

The classes of proportional elements (i.e. rays) are the equivalence classes of
these preorder relations. The importance of these relations for the geometry of
max cones is expressed by the following principle. Denote $\supp(y):=\{i\mid y_i\neq 0\}$.

\begin{proposition}
\label{p:multiorder}
Let $V=\spann(S)$, $S\subseteq\Rpn$. Then the following
are equivalent.
\begin{itemize}
\item[1.] $y\in V$.
\item[2.] For all $j\in\supp(y)$ there exists $v\in S$ such that $v\leq_j y$.
\end{itemize}
\end{proposition}

This principle appeared as a set covering condition, see
Proposition~ \ref{setcov} below, already in the works of Vorobyev
\cite{Vor-67} and Zimmermann \cite{Zim-76}, and in the above form
(or with a subtle difference) it appeared quite recently in the
works of Joswig \cite{Jos-05}, Ni\c{t}ic\u{a} and Singer
\cite{NS-07b}, and Butkovi\v{c} et al. \cite{BSS-07}, see also
\cite{But-03} and \cite{DS-04}.

As it was remarked by Ni\c{t}ic\u{a} and Singer \cite{NS-07b}, the above proposition
means that the geometry of max cones is a special case of the {\em multiorder convexity}
\cite{MLS-91}. In the multiorder convexity, one has a set of order relations,
and a point $y$ is said to belong to the convex hull of $S$, if for any order
there is a point in $S$ which precedes $y$ with respect to that order.

The following proposition is the Tropical Carath\'{e}odory Theorem,
see Helbig \cite{Hel-88}, Develin and Sturmfels \cite{DS-04}, and
also \cite{BSS-07,GM-08}. Note that it follows from Proposition
\ref{p:multiorder}.
\begin{proposition}
\label{cara}
Let $S\subseteq\Rpn$. Then $y\in\spann(S)$ if and only if there exist $k$ vectors
$v^1,\ldots, v^k\in S$, where $k=|\supp(y)|$, such that $y\in\spann(v^1,\ldots,v^k)$.
\end{proposition}

The multiorder principle also means the following description of extremals \cite{BSS-07}.

\begin{proposition}
\label{ex-minima}
Let $V\subseteq\Rpn$ be a max cone generated by $S$ and let $v\in V$, $v\neq 0$.
Then the following are equivalent.
\begin{itemize}
\item[1.] $v$ is an extremal in $V$.
\item[2.] For some $j\in\supp(v)$, $v$ is minimal with respect to $\leq_j$
in $V$.
\item[3.] For some $j\in\supp(v)$, $v$ is minimal with respect to $\leq_j$
in $S$.
\end{itemize}
\end{proposition}

Propositions \ref{p:multiorder} and \ref{ex-minima} lead to a number
of statements about generators, extremals and bases of max cones
\cite{BSS-07}, we mention only the following two of them. An element
$u\in\Rpn$ is called {\em scaled}, if $||u||=1$, where $||\cdot||$
denotes some fixed norm (say, the ordinary norm or the max norm).
For the following proposition see Butkovi\v{c} et al. \cite{BSS-07},
and also \cite{DS-04,Wag-91} for closely related
statements.

\begin{proposition}
\label{basis-unique}
Let $E$ be the set of scaled extremals in a max cone $V\subseteq\Rpn$ and let
$S\subseteq\Rpn$ consist of scaled elements. Then the following are equivalent.
\begin{itemize}
\item[1.] The set $S$ generates $V$ and none of the elements in $S$ are redundant.
\item[2.] $S=E$ and $S$ generates $V$.
\item[3.] The set $S$ is a basis for $V$.
\end{itemize}
\end{proposition}

Proposition \ref{basis-unique} means that if a scaled basis of a max
cone exists, then it is unique and consists of all scaled extremals,
i.e., all the elements that are minimal with respect to some
preorder relation $\leq_i$. In particular, a scaled basis of a
finitely generated max cone $V$ exists and is unique, and the
cardinality of this basis will be called the {\em max-algebraic
dimension} of $V$.

The following result is analogous to Minkowski's theorem about
extremal points of convex sets, and was obtained independently by
Gaubert and Katz \cite{GK-07} and Butkovi\v{c} et al. \cite{BSS-07}.

\begin{proposition}
\label{mink} Let $V\subseteq\Rpn$ be a closed max cone. Then $V$ is
generated by its set of extremals, and any vector in $V$ is a max
combination of no more than $n$ extremals.
\end{proposition}

Note that any finitely generated max cone is closed (\cite{BSS-07,Jos-05}).
One may also think of colourful extensions of Propositions \ref{cara} and \ref{mink}
in the sense of B\'{a}r\'{a}ny \cite{Bar-82}, and progress in this direction is due
to Gaubert and Meunier \cite{GM-08}, see also Theorem \ref{colour-mink} below.

\subsection{Projectors and separation}

Given a closed max cone $V\subseteq\Rpn$, we can define a {\it nonlinear projector} $P_V$ by
\begin{equation}
\label{proj-def}
P_V(y)=\max\{v\in V\mid v\leq y\}.
\end{equation}
This operator is homogeneous: $P_V(\lambda y)=\lambda P_V(y)$,
isotone: $y^1\leq y^2\Rightarrow P_V(y^1)\leq P_V(y^2)$,
nonincreasing: $P_V(y)\leq y$, and continuous, see \cite{CGQS-05}
for the proof. For any vector $y$ there are coordinates which do not
change under the action of the projector: $P_V(y)_i=y_i$. These
coordinates will be called {\em sleepers}. Projectors lead to
separation theorems of the following kind, see
\cite{CGQS-05,DS-04,GS-07,Jos-05} and introduction for some
historical remarks.

\begin{theorem}
\label{separ0}
Let $V\subseteq\Rpn$ be a closed max cone and
let $y\in\Rpn$ be not in $V$. Then there exist a positive vector $\Tilde{y}$
and a max cone $\Tilde{V}\supseteq V$ containing positive vectors such that the set
\begin{equation}
\label{halfs0}
\halfs=\{v\mid\bigoplus_{i=1}^n \Tilde{y}_i^{-1} v_i\geq\bigoplus_{i=1}^n
(P_{\Tilde{V}}(\Tilde{y}))^{-1}_i v_i\}
\end{equation}
contains $V$ but not $y$. If $y$ is positive and $V$ contains positive vectors, then one can take
$\Tilde{y}=y$ and $\Tilde{V}=V$.
\end{theorem}

The set $\halfs$ defined in \eqref{halfs0} is an instance of the max analogue of a halfspace,
which is generally a set of the form
%\begin{equation*}
%\label{halfs-gen}
$\{v\mid\bigoplus_{i=1}^n u^1_iv_i\geq \bigoplus_{i=1}^n u^2_iv_i\}.$
%\end{equation*}

By comparing this to \eqref{halfs0} we see that a separating halfspace has both $u^1$ and $u^2$ positive
and $u^1\leq u^2$, so that the inequality in \eqref{halfs0} can be replaced by equality:

\begin{equation}
\label{halfs01}
\halfs=\{v\mid\bigoplus_{i=1}^n \Tilde{y}_i^{-1} v_i=
\bigoplus_{i=1}^n (P_{\Tilde{V}}(\Tilde{y}))^{-1}_i v_i\}.
\end{equation}

The relation of Theorem \ref{separ0} to the multiorder principle was made explicit by Joswig \cite{Jos-05}.
Denote, for any positive $y$,
%\begin{equation}
%\label{Delta}
$\Delta_i(y)=\{u\in\Rpn\mid u\leq_i y\}.$
%\end{equation}
Observe that $\bigcup_{i=1}^n \Delta_i(y)=\Rpn$, and that the separating
halfspace defined by \eqref{halfs0} or equivalently \eqref{halfs01}
can also be written as
\begin{equation}
\label{sectors}
\halfs=\bigcup_{i\in\sle(P_{\Tilde{V}},\Tilde{y})} \Delta_i(P_{\Tilde{V}}(\Tilde{y})),
\end{equation}
where $\sle(P_{\Tilde{V}},\Tilde{y})$ is the set of sleepers, i.e., the indices $k$
such that $(P_{\Tilde{V}}(\Tilde{y}))_k=\Tilde{y}_k$.
Thus, in terms of the multiorder, the separation theorem says that,
given a point $y$ and a closed max cone $V$, there is a point $P_{\Tilde{V}}(\Tilde{y})$
such that the union of some sectors $\Delta_i(P_{\Tilde{V}}(\Tilde{y}))$ contains the
whole $V$ while the complement of this union contains $y$.

If a max cone is generated by the columns of a matrix $A\subseteq\Rpnm$,
then, denoting $P_A:=P_{\spann(A)}$, we deduce from \eqref{proj-def} that
\begin{equation}
\label{CG-proj}
P_A(y)=A\otimes(\overline{A}\otimes' y),
\end{equation}
where $\overline{A}$ is the {\em Cuninghame-Green inverse} of $A$ defined by
$\overline{a}_{ij}=a_{ji}^{-1}$, and $\otimes'$ denotes the min-times
matrix product. When calculating \eqref{CG-proj}, we put by
convention that $0^{-1}=\infty$ and $0\otimes +\infty=0$. In this
form \eqref{CG-proj}, the nonlinear projectors were studied by
Cuninghame-Green \cite{CG:79}. We also note that formula \eqref{CG-proj} 
represents a projector as a min-max function 
in the sense of \cite{CZ-02,CTGG-99,OP-97}, with addition being replaced by
multiplication.

When $V$ is an arbitrary closed max cone, $P_V$ can be expanded in infinite sum of 'elementary'
projectors using the following 'scalar product', or an instance of residuation \cite{CGQ-04,CGQS-05}:
\begin{equation*}
%\label{y/x}
y/v:=\min_{i\in\supp(v)} y_i v_i^{-1}=\max\{\lambda\mid \lambda
v\leq y\}.
\end{equation*}
Namely,
\begin{equation}
\label{pv(x)}
P_V(y)=\bigoplus_{v\in V} y/v\ v.
\end{equation}
Formula \eqref{CG-proj} is a special case of \eqref{pv(x)}, when $V$
is finitely generated. Using the multiorder, we can obtain the
following refinement of \eqref{pv(x)}. Denote by $\wedge$ the
componentwise minimum of vectors in $\Rpn$.

\begin{theorem}
\label{projector1} Suppose that $V\subseteq\Rpn$ is a closed max
cone. Then for any $y\in\Rpn$, the components $(P_V(y))_i$, for
$i\in\supp(y)$, are equal to
\begin{equation}
\label{pv(y)i}
(P_V(y))_i=\bigoplus_{v\in E_i} y/v\ v_i,
\end{equation}
where $E_i$ is the set of scaled points of $V$, minimal with respect to $\leq_i$.
The projector $P_V$ is linear with respect to the componentwise minimum
$\wedge$ if and only if every set $E_i$ is a singleton.
\end{theorem}

\begin{proof}
Writing \eqref{pv(x)} componentwise, we have that
\begin{equation*}
\label{pv(y)i2}
(P_V(y))_i=\max_{v\in V:\, v_i\neq 0}(v_i\min_{k: v_k\neq 0} y_k v_k^{-1})=
\max_{v\in V:\, v_i\neq 0}\min_{k: v_k\neq 0} y_k(v_k v_i^{-1})^{-1}.
\end{equation*}
By Proposition \ref{mink}, any closed max cone has a scaled basis
$E$. Denote by $E_i$ the set of scaled vectors minimal with respect
to $\leq_i$, then for all $v\in V$ and any $i\in\supp(v)$ there is
$v^i\in E_i$ such that $v^i\leq_i v$ and hence $(v_k^i
(v_i^i)^{-1})^{-1}\geq (v_k v_i^{-1})^{-1}$ for all $k$. This proves
\eqref{pv(y)i}, and \eqref{pv(y)i} implies that if all the sets
$E_i$ consist of one element, then the projector is expressed by a
min-times matrix. Now suppose that there is an $i$ such that $E_i$
has at least two elements, say, $u$ and $v$. Then $P_V(u)=u$ and
$P_V(v)=v$. If the projector is linear with respect to the
componentwise minimum $\wedge$, then $P_V(uu_i^{-1}\wedge vv_i^{-1})
=uu_i^{-1}\wedge vv_i^{-1}$, hence $w=uu_i^{-1}\wedge vv_i^{-1}\in
V$. As $w_i=1$, we have that $w\leq_i v$ and $w\leq_i u$. As $u$ and
$v$ are both minimal with respect to $\leq_i$, $w$ is not equal to
either of them, which leads to a contradiction with the minimality
of $u$ and $v$. The proof is complete.
\end{proof}

\section{The role of Kleene stars}

\subsection{Kleene stars and Kleene cones}

We start this section with some necessary definitions.
Let $A=(a_{ij})\in\Rpnn$. The weighted digraph $D_A=(N(A),E(A))$,
whose nodes are $N(A)=[n]$ and whose edges $E(A)=N(A)\times N(A)$
have weights $w(i,j)=a_{ij}$, is called the {\em digraph
associated} with $A$.
Suppose that $\pi =(i_{1},...,i_{p})$ is a path in $D_A$, then the \textit{%
weight} of $\pi $ is defined to be $w(\pi
,A)=a_{i_{1}i_{2}}a_{i_{2}i_{3}}\ldots a_{i_{p-1}i_{p}}$ if $p>1$, and $0$ if $p=1$.
A path which begins at $i$ and ends at $j$ will be called an {\em $i\to j$ path}.
If the starting node of a path coincides with the end node then the path
is called a {\it cycle}.

A path $\pi$ is called {\em positive} if $w(\pi,A)>0$. If for all $i,j\in[n]$ there
exists a positive $i\to j$ path, then $A$ is called {\em irreducible}.

The \textit{maximum cycle geometric mean} of $A$, further denoted by $\lambda(A)$,
is defined by the formula
\begin{equation*}
\lambda (A)=\max_{\sigma }\mu (\sigma ,A),  \label{mcm}
\end{equation*}%
where the maximisation is taken over all cycles in the digraph and
\begin{equation*}
\mu (\sigma ,A)=w(\sigma ,A)^{1/k}  \label{cm}
\end{equation*}%
denotes the \textit{geometric mean} of the cycle $\sigma =(i_{1},...,i_{k},i_{1})$.

The following fact was proved by Carr\'{e} \cite{Car-71}, see also
\cite{BCOQ,CG:79}.

\begin{proposition}
\label{Carre}
Let $A\in \Rpnn$. The series
\begin{equation}
\label{e:closeries}
A^*=I\oplus A\oplus A^2\oplus\ldots
\end{equation}
converges to a finite limit and is equal to $I\oplus A\oplus\ldots\oplus A^{n-1}$ if and only
if $\lambda(A)\leq 1$. In this case also $\lambda(A^*)\leq 1$.
\end{proposition}

The matrix series $A^*$ defined by \eqref{e:closeries} is called the
{\em Kleene star} of $A$, which comes from the theory of automata, see
Conway \cite{Con:71}. Kleene stars enjoy the property $(A^*)^2=A^*$, i.e., they are
{\em multiplicatively idempotent}. Their diagonal entries are all equal to $1$, i.e., the
Kleene stars are {\em increasing}. Actually these two properties are also sufficient for a matrix to
be a Kleene star, and further by a Kleene star we will also mean any matrix
with these two properties. We also note that $(A^*)^2=A^*$ implies that $(A^*)^*=A^*$.

A max cone will be called a {\em Kleene cone} if it can be represented as
max-algebraic column span of a Kleene star.

In terms of the multiorder, we can say that a matrix $A$ is a Kleene
star if and only if $a_{ii}=1$ for all $i\in[n]$ and $A_{\cdot
i}\leq_i A_{\cdot k}$ for all $i,k$ such that $a_{ik}\neq 0$. That
is, $A$ is a Kleene star if and only if $a_{ii}=1$ and $A_{\cdot i}$
is the unique minimum of $\spann(A^*)$ with respect to $\leq_i$ for
all $i\in[n]$, so that all the sets $E_i$ defined in Theorem~\ref{projector1} are singletons.
The last sentence of Theorem~\ref{projector1} can be now formulated as follows.

\begin{proposition}
\label{proj-linear} $P_V$ is a min-times linear operator if
and only if $V$ is a Kleene cone. If $V=\spann(A)$, where $A$ is a
Kleene star, then $P_V(y)=\overline{A}\otimes' y$ for all $y$.
\end{proposition}

Kleene stars play crucial role in the description of max-algebraic
eigenvectors and subeigenvectors of nonnegative matrices. If for
some $x$ and $\lambda$ we have that $A\otimes x=\lambda x$, then
$\lambda$ is a {\em max-algebraic eigenvalue} of $A$, and $x$ is a
{\em max-algebraic eigenvector} associated with this eigenvalue.
Analogously, $x$ is called a {\em max-algebraic subeigenvector} associated with
$\lambda$, if $A\otimes x\leq\lambda x$.

The well-known Perron-Frobenius theorem has a max-algebraic analogue \cite{BCOQ,Bap-98,CG:79,Vor-67}.
\begin{theorem}
\label{PF-maxalg} Let $A\in\Rpnn$.
\begin{itemize}
\item[1.] $A$ has a max-algebraic eigenvalue, and the number of such eigenvalues is less
than or equal to $n$.
%\item[2.] $A$ has a positive eigenvector associated with the largest eigenvalue, and
%any positive eigenvector is associated with the largest eigenvalue;
\item[2.] $\lambda(A)$ is the largest eigenvalue of $A$.
\item[3.] If $A$ is irreducible, then $\lambda(A)$ is the unique max-algebraic eigenvalue of $A$
and all eigenvectors associated with $\lambda(A)$ are positive.
\end{itemize}
\end{theorem}

The set of eigenvectors associated with a fixed
eigenvalue $\lambda$ is a max cone, and analogously the set of subeigenvectors
associated with a fixed $\lambda$ is a max cone, so they will be called the {\em eigencone}
and the {\em subeigencone} associated with $\lambda$.
For a nonnegative square matrix $A\in\Rpnn$ the eigencone associated with $1$
will be denoted by $V(A)$, and the subeigencone associated with $1$ will be
denoted by $V^*(A)$. A matrix $A\in\Rpnn$ is called {\em definite}, if $\lambda(A)=1$.
We do not lose much generality when considering definite matrices, as for any matrix $A$
with $\lambda(A)\neq 0$, the matrix $A/\lambda(A)$ is definite and has the same
eigenvectors and subeigenvectors as $A$.

Any subeigencone is a Kleene cone, and the other way around.
\begin{proposition}
\label{subeig-kleene}
Let $A\in\Rpnn$ be definite, then $V^*(A)=V(A^*)=\spann(A^*)$.
\end{proposition}
\begin{proof}
First note that by Proposition~\ref{Carre}, if $\lambda(A)=1$ then
$A^*$ exists and $\lambda(A^*)=1$.

We show that $V^*(A)=V(A^*)$. Suppose that $A^*\otimes x=x$,
then $A\otimes x\leq x$, because $A\leq A^*$. If $A\otimes x\leq x$, then $(I\oplus A)\otimes x=x$
and also $A^*\otimes x=x$, since $A^m\otimes x\leq x$ for any $m$
(due to the isotonicity of matrix multiplication).

We show that $V(A^*)=\spann(A^*)$. It is immediate that
$V(A^*)\subseteq\spann(A^*)$, as $V(A)\subseteq\spann(A)$ for any
matrix $A$. If $A^*$ converges, then $A\otimes A^*=A\oplus
A^2\oplus$..., so $A\otimes A^*\leq A^*$ meaning that each column of
$A^*$ is a subeigenvector of $A$. Hence $\spann(A^*)\subseteq
V^*(A)$.
\end{proof}

The positivity of subeigenvectors is addressed in the following observation.

\begin{proposition}
\label{subeig-pos}
Let $A\subseteq\Rpnn$ be such that $a_{ii}=1$ for all $i\in[n]$.
Then $V^*(A)$ contains a positive vector if and only if $A$ is definite.
\end{proposition}
\begin{proof}
The ``if'' part: If $A$ is definite, then by
Proposition~\ref{subeig-kleene} $V^*(A)=\spann(A^*)$ and we can
take, for a positive subeigenvector of $A$, any max combination of
all the columns of $A^*$ with positive coefficients.

The ``only if'' part: Suppose that there exists a positive $x$ such
that $A\otimes x\leq x$, and take a cyclic permutation
$\tau=(i_1,\ldots,i_k)$ of a subset of $[n]$. Then we have that
$a_{i_li_{l+1}}x_{i_{l+1}}\leq x_{i_l}$ for $l\in[k]$, assuming
$i_{k+1}:=i_1$. Multiplying all these inequalities and cancelling the coordinates of $x$ we have that
$w(\tau,A)\leq 1$. Hence $\lambda(A)\leq 1$. As all diagonal entries
are equal to $1$, we have that $\lambda(A)=1$.
\end{proof}

Proposition \ref{subeig-kleene} implies that if $A$ is a Kleene
star, then
\begin{equation*}
\label{e:afriat1} \spann(A)=V(A)=V(A^*)=V^*(A)=\{x\mid a_{ij} x_j\leq x_i,\
i,j\in[n]\},
\end{equation*}
and it is not hard to see the following.

\begin{proposition}
\label{p:afriat2} Let $K$ be a max cone in $\Rpn$. Then it is a
Kleene cone if and only if for some matrix $B$ it is the solution
set of the system of inequalities $b_{ij} x_j\leq x_i,\ i,j\in[n],$
%\begin{equation*}
%\label{e:afriat2}
%K=\{x\mid b_{ij} x_j\leq x_i,\ i,j\in[n]\},
%\end{equation*}
satisfied by at least one positive $x$.
\end{proposition}
\begin{proof}
The ``if'' part: If the system is satisfied by a positive $x$, then
$b_{ii}\leq 1$ for all $i\in[n]$. Take $\Tilde{B}:=I\oplus B$, then
$\Tilde{B}$ has all diagonal entries equal to $1$,
$K=V^*(\Tilde{B})$ and there is a positive $x\in V^*(\Tilde{B})$. By
Proposition~\ref{subeig-pos}, $\Tilde{B}$ is definite, and by
Proposition~\ref{subeig-kleene}, $K=\spann((\Tilde{B})^*)$.

The ``only if'' part: If $K$ is a Kleene cone $\spann(A^*)$, then by
Proposition~\ref{subeig-kleene} and Proposition~\ref{subeig-pos} we
can take $B:=A^*$.
\end{proof}

%System \eqref{e:afriat2} is a multiplicative version of an Afriat-Varian system
%of inequalities \cite{Roc:70},
The above observations imply that Kleene cones
are convex cones, and that they
have many close relatives in the realm of combinatorial geometry, see
Joswig and Kulas \cite{JK-08}.

One may think of various systems of inequalities describing the same Kleene cone. However, the
Kleene star which defines this cone is unique \cite{Ser-07}.

\begin{proposition}
\label{clounique}
Suppose that $A$ and $B$ are two Kleene stars. Then $A=B$ if and only if $\spann(A)=\spann(B)$.
\end{proposition}

We now describe the bases of $V(A)$ and $V^*(A)$, for a definite matrix $A\in\Rpnn$.
The cycles with the cycle geometric mean equal to $1$ are
called {\em critical}, and the nodes and the edges of $D_A$ that
belong to critical cycles are called {\em critical}. The set of
critical nodes is denoted by $N_c(A)$, the set of critical
edges is denoted by $E_c(A)$, and the {\em critical digraph} of
$A$, further denoted by $C(A)=(N_c(A),E_c(A))$, is the digraph that
consists of all critical nodes and critical edges of $D_A$. All
cycles of $C(A)$ are critical \cite{BCOQ}.
For two vectors $x$ and $y$, we write $x\sim y$ if
$x=\lambda y$ for $\lambda>0$. The following theorem follows from
well-known results on the max-algebraic spectral theory \cite{BCOQ,CG:79,GauThesis}.

\begin{theorem}
\label{subeigs}
Let $A\in\Rpnn$ be definite, and let
$M(A)$ denote a set of indices such that for
each strongly connected component of $C(A)$ there is a unique index
in $M(A)$ which belongs to that component.
\begin{itemize}
\item[1.] The following statements are equivalent: $A^*_{\cdot i}\sim A^*_{\cdot j}$,
$A^*_{i\cdot}\sim A^*_{j\cdot}$, $i$ and $j$ belong to the same strongly connected component
of $C(A)$.
\item[2.] Any column of $A^*$ is a max extremal of $\spann(A^*)$.
\item[3.] The subeigencone of $A$, which is the eigencone of $A^*$, is
\begin{equation*}
V^*(A)=V(A^*) =\left\{\bigoplus_{i\in M(A)} \alpha_i A^*_{\cdot i}\oplus\bigoplus_{j\notin
C(A)} \alpha_j A^*_{\cdot  j},\ \alpha_i,\alpha_j\in\Rp\right\},
\end{equation*}
and none of the columns of $A^*$ in this description are redundant.
\item[4.] The eigencone of $A$ is
\begin{equation*}
V(A) =\left\{\bigoplus_{i\in M(A)} \alpha_i A^*_{\cdot  i},\
\alpha_i\in\Rp\right\},
\end{equation*}
and none of the columns of $A^*$ in this description are redundant.
\end{itemize}
\end{theorem}

Proposition \ref{basis-unique} and Theorem \ref{subeigs} imply that
extremals of $V^*(A)$ are precisely the columns of $A^*$, so the
columns of $A^*$, after eliminating the proportional ones,
constitute the basis of $V^*(A)=\spann(A^*)$, and the columns whose
indices belong to $C(A)$ constitute the basis of $V(A)$.
Denote by $n_c(A)$ the number of strongly connected
components in $C(A)$, and denote by $\overline{N_c(A)}$ the set of nodes that are
not critical. Theorem \ref{subeigs} yields the following
corollary.

\begin{proposition}
\label{maxdim}
For any definite matrix $A\in\Rpnn$, the max-algebraic dimension
of the subeigencone of $A$ is equal to $n_c(A)+|\overline{N_c(A)}|$. The max-algebraic dimension
of the eigencone is equal to $n_c(A)$.
\end{proposition}

Kleene cones are both convex cones and max cones. They are inhabitants of two worlds,
that of max algebra and tropical convexity, and that of nonnegative linear algebra and ordinary
convexity. One might think of an interplay between these worlds. For a definite
matrix $A$, define the linear space
\begin{equation}
\Lin(C(A))=\{x\in\Rn\mid a_{ij} x_j= x_i,\ (i,j)\in E_c(A)\}.
\end{equation}
A proof of the following theorem can be found in \cite{SSB-08}.
\begin{theorem}
\label{subeig-dimension}
Let $A\in\Rpnn$ be a definite matrix. Then $\Lin(C(A))$ is the linear hull of the convex cone
$V^*(A)$. The linear dimension of $V^*(A)$, i.e., the dimension of $\Lin(C(A))$, is equal to the
max-algebraic dimension of $V^*(A)$, i.e., to $n_c(A)+|\overline{N_c(A)}|$.
\end{theorem}

The intersection of Kleene cones is again a Kleene cone. More
precisely, we have the following proposition, see Butkovi\v{c}
\cite{But-85} for the case $k=2$. The proof is based
on the formula $(A^*\oplus B^*)^*=(A^*\otimes B^*)^*$, which follows from
$(A\oplus B)^*=A^*\otimes(B\otimes A^*)^*$
\cite{Con:71}, and on the observations above.

\begin{proposition}
\label{star-inters-k} Let $A^{(1)},\ldots,A^{(k)}\in\Rpnn$
be Kleene stars. The following are equivalent.
\begin{itemize}
\item[1.] $\bigcap_{i=1}^k \spann(A^{(i)})$ contains a positive vector.
\item[2.] $\lambda(\bigoplus_{i=1}^k A^{(i)})=1$.
\item[3.] $\lambda(\bigotimes_{i=1}^k A^{(\pi(i))})=1$ for some permutation
$\pi$ of $\{1,\ldots,k\}$.
\item[4.] $\lambda(\bigotimes_{i=1}^k A^{(\pi(i))})=1$ for all permutations
$\pi$ of $\{1,\ldots,k\}$.
\end{itemize}
If any of these equivalent conditions are true, then
\begin{equation}
\label{e:star-inters}
\bigcap_{i=1}^k\spann(A^{(i)})=\spann((\bigoplus_{i=1}^k A^{(i)})^*)=
\spann((\bigotimes_{i=1}^k A^{(\pi(i))})^*)
\end{equation}
for all permutations $\pi$.
\end{proposition}
\begin{proof}
Complete $\R_{\max,\times}$ with $+\infty$ and assume $a\times
+\infty=+\infty$ for any positive $a$ and $0\times +\infty=0$.
Matrix algebra over this completed semiring is a regular algebra in
the sense of \cite{Con:71}. This means in particular that $A^*$ is
always defined, $(A^*)^*=A^*$, $(A\oplus B)^*=A^*\otimes(B\otimes
A^*)^*$ and $(A\otimes B)^*=I\oplus (A\otimes(B\otimes A)^*)$. If
$A$ and $B$ are two Kleene stars, then
\begin{equation}
\label{AB*}
\begin{split}
(A\otimes B)^* & =I\oplus (A\otimes(B\otimes A)^*)=A\otimes(B\otimes
A)^*=\\ &=(A\oplus B)^*=(B\oplus A)^*=(B\otimes A)^*.
\end{split}
\end{equation}
It can be shown by induction that $(A^{(1)}\oplus\ldots\oplus
A^{(k)})^*=(A^{(\pi(1)}\otimes\ldots\otimes A^{(\pi(k)})^*$ for any
permutation $\pi$ of $\{1,\ldots,k\}$. Using
Proposition~\ref{Carre} we obtain that $\lambda(\bigoplus_{i=1}^k A^{(i)})\leq 1$
is true if and only if $\lambda(\bigotimes_{i=1}^k A^{\pi(i)})\leq
1$ is true for some $\pi$, and hence if and only if the same is true
for all $\pi$. The inequalities here can be replaced by equalities,
since all diagonal entries, and hence all eigenvalues, of any
product or entrywise maximum of Kleene stars, are greater than or
equal to $1$. This yields equivalence of 2., 3., and 4.

We now prove the equivalence between 1. and 2., and
\eqref{e:star-inters}. We have that
\begin{equation}
\label{star-inters-imm}
%V(\bigoplus_{i=1}^k A^{(i)})=
V^*(\bigoplus_{i=1}^k A^{(i)})=\bigcap_{i=1}^k V^*(A^{(i)})=
\bigcap_{i=1}^k \spann(A^{(i)}),
\end{equation}
where the first equality is immediate, and the second equality
follows from Proposition~ \ref{subeig-kleene}. Note that all
diagonal entries of $\bigoplus_{i=1}^k A^{(i)}$ are $1$, and by
Proposition~ \ref{subeig-pos}, $V^*(\bigoplus_{i=1}^k A^{(i)})$
contains a positive vector if and only if $\lambda(\bigoplus_{i=1}^k
A^{(i)})=1$. This, together with \eqref{star-inters-imm}, implies
the equivalence between assertions 1. and 2. By
Proposition~\ref{subeig-kleene}, $V^*(\bigoplus_{i=1}^k
A^{(i)})=\spann((\bigoplus_{i=1}^k A^{(i)})^*)$ since
$\lambda(\bigoplus_{i=1}^k A^{(i)})=1$, which yields~
\eqref{e:star-inters}.
\end{proof}

\subsection{Cellular decomposition}

We have described some properties of Kleene cones. Though such cones
are very special, they can be viewed as building blocks, or atoms,
of any finitely generated max cone. This can be seen as the main
idea of the cellular decomposition, an ingenuous concept of Develin
and Sturmfels \cite{DS-04}, which we adjust below to the setting of
max cones.

Let $A\subseteq\R_+^{n\times m}$ be a nonnegative matrix with $m$
nonzero columns and $n$ nonzero rows. The {\em column type} of $y$
with respect to $A$ is defined to be the $m$-tuple of subsets
$T_1,\ldots, T_m$ of $[n]$, where every $T_j$, for $j\in[m]$ is
defined by
\begin{equation*}
\label{column type} T_j=\{i\in [n]\mid a_{ij} y_i^{-1}\geq a_{kj}
y_k^{-1},\ k\in [n]\}=\{i\in[n]\mid y\geq_i A_{\cdot j}\}.
\end{equation*}
The {\em row type} of $y$ with respect to $A$ is an $n$-tuple of
subsets $S_1,\ldots,S_n$ of $[m]$, where every $S_i$, for $i\in[n]$,
is defined by
\begin{equation*}
\label{row type}
\begin{split}
S_i&=\{j\in [m]\mid a_{ij} y_i^{-1}\geq a_{kj} y_k^{-1},\ k\in
[n]\}=\{j\in[m]\mid y\geq_i A_{\cdot j}\}=\\ &=\{j\in [m]\mid i\in
T_j\}.
\end{split}
\end{equation*}
The theory of $A\otimes x=y$ systems
\cite{BCOQ,But-03,CG:79,DS-04,Vor-67,Zim-76} is based on the
following {\em set covering conditions} for $y$ to be in
$\spann(A)$, see the proposition below. The multiorder principle
(Proposition~\ref{p:multiorder}) can be seen as a reformulation of
these conditions, therefore we leave the proposition below without proof.
\begin{proposition}
\label{setcov} Let $A\in\Rpnm$ have all rows and columns nonzero and
let $y\in\Rpn$ be a positive vector with the column type
$T=(T_1,\ldots,T_m)$ and the row type $S=(S_1,\ldots,S_n)$. The
following are equivalent.
\begin{itemize}
\item[1.] $y\in\spann(A)$;
\item[2.] $\bigcup_{i=1}^m T_i=[n]$;
\item[3.] none of $S_i, i\in [n]$ are empty.
\end{itemize}
\end{proposition}

See also Akian et al. \cite{AGK-04} for an
infinite-dimensional generalisation in the context of Galois
connections.

Following Develin and Sturmfels \cite{DS-04}, we can see this from a
geometric viewpoint. For any row type $S$, we define its {\em region} with
respect to $A$ by
\begin{equation*}
\label{xsdef} X_S=\{ y\ \text{positive}\mid y_k y_i^{-1}\geq a_{kj}
a_{ij}^{-1},\ \forall k,i,\ \forall j\in S_i\}.
\end{equation*}
Proposition \ref{setcov} means that the part of $\spann(A)$ consisting of all
positive vectors is the union of
the regions $X_S$ such that $S$ do not contain empty sets (\cite{DS-04}, Theorem 15).
If $X_S$ is not empty, then the closure of $X_S$ is
\begin{equation}
\label{cloxsdef} \operatorname{cl}(X_S)=\{y\in\Rpn\mid a_{kj}
a_{ij}^{-1} y_i\leq y_k,\ \forall k,i,\ \forall j\in S_i\}.
\end{equation}
It follows from the results of \cite{DS-04} that the relative
interiors of regions build up a cellular decomposition of the
positive part of $\Rpn$. We will need a weaker statement, but
without positivity.
\begin{proposition}
\label{cellular1} Suppose that $A\in\Rpnm$ has all rows and columns
nonzero. Then the max cone $\spann(A)$ is the union of
$\operatorname{cl}(X_S)$ such that $X_S$ are not empty and $S$ do
not contain empty sets.
\end{proposition}
\begin{proof}
As $A$ has all rows nonzero, the max cone $\spann(A)$ contains positive vectors. By Proposition
\ref{setcov} if $y$ is positive, then $y\in\spann(A)$ if and only if
the row type of $y$ does not contain empty sets. Hence the positive
part of $\spann(A)$ is the union of nonempty $X_S$ such that $S$ do
not contain empty sets. Further, $\spann(A)$ is the closure of its
positive part. Indeed, $\spann(A)$ contains positive vectors and for
any $u\in\spann(A)$ and a positive $v\in\spann(A)$ we can take
$w=u\oplus\varepsilon v\in\spann(A)$, so that
$||w-u||\leq\varepsilon ||v||$ (the max norm) and $w$ is positive.
Hence $\spann(A)$ is the union of closed regions
$\operatorname{cl}(X_S)$ such that $X_S$ are not empty and $S$ do
not contain empty sets.
\end{proof}

From the max-algebraic point of view, an important role in the cellular decomposition
is played by {\em strongly definite matrices}, which are definite matrices
with all diagonal entries equal to $1$. Note that any Kleene star is a strongly definite matrix.

Observe that $\operatorname{cl}(X_S)$ is the subeigencone of the $n\times n$
matrix $A^S=(a^S_{ij})$ defined by
\begin{equation}
\label{asdef}
a^S_{ij}=
\begin{cases}
{\bigoplus_{k\in S_j}} a_{ik} a_{jk}^{-1}, & \text{if $S_j\ne\emptyset$},\\
\kroneker_{ij}, & \text{if $S_j=\emptyset$,}
\end{cases}
\end{equation}
where $\kroneker_{ij}$ are Kronecker symbols ($\kroneker_{ij}=0$ if $i\neq j$ and $\kroneker_{ij}=1$
if $i=j$). It is immediate that all diagonal entries
of $A^S$ are equal to $1$.
We have the following proposition which can be used to compute the
generators of any closed region, a preliminary version of this proposition
appeared in \cite{Ser-07}.

\begin{proposition}
\label{regionskleene} The closed region $\operatorname{cl}(X_S)$
contains positive vectors if and only if $A^S$ is a strongly
definite matrix, and in this case
$\operatorname{cl}(X_S)=V^*(A^S)=\spann((A^S)^*)$.
\end{proposition}
\begin{proof} From \eqref{cloxsdef} and \eqref{asdef} one infers
that $\operatorname{cl}(X_S)=V^*(A^S)$. After that, the claim
follows from Proposition~\ref{subeig-kleene} and
Proposition~\ref{subeig-pos}.
\end{proof}

Propositions \ref{cellular1} and \ref{regionskleene} have the following
consequences.

\begin{proposition}
\label{cellular2} For any matrix $A\in\Rpnm$ with no zero rows there
exist Kleene stars $A^{(1)},\ldots,A^{(l)}\in\Rpnn$ such that
%\begin{equation}
$\spann(A)=\bigcup_{i=1}^l \spann(A^{(i)}).$
%\end{equation}
\end{proposition}

\begin{proposition}
\label{cellular3} For any matrix $A\in\Rpnm$ with no zero rows there
exist Kleene stars $A^{(1)},\ldots, A^{(l)}\in\Rpnn$ such that for
any $y\in\Rpn$ we have that $P_A y=\overline{A^{(k)}}\otimes' y$ for
some $k$.
\end{proposition}

To express the dimension of a region, Develin and Sturmfels
\cite{DS-04} introduce the undirected graph $G_S$: The set of nodes
of this graph is $[n]$, it contains all loops $(i,i)$, and for $i\neq
j$ an edge $(i,j)$ belongs to $G_S$ if and only if there exists
$k\in S_i\cap S_j$. The following observation relates
this notion to max algebra.

\begin{proposition}
\label{gs=gc} Let $A\in\Rpnm$ be a matrix with no zero rows and
columns, let $y\in\Rpn$ be a positive vector and $S$ be the row type
of $y$ with respect to $A$. Then $G_S=C(A^S)$.
\end{proposition}
\begin{proof}
Note that as all entries of $A^S$ are equal to $1$, the graph $C(A^S)$ contains all loops.

Let $i\neq j$ and $(i,j)\in G_S$, then there exists $k\in S_i\cap
S_j$. It follows that $a_{ik} a_{jk}^{-1}=y_i y_j^{-1}\geq a_{il}
a_{jl}^{-1}$ for all $l\in S_j$, and therefore
$a_{ij}^S=a_{ik}a_{jk}^{-1}$. Analogously,
$a_{ji}^S=a_{jk}a_{ik}^{-1}$, and therefore $a_{ij}^S a_{ji}^S=1$ so
that $(i,j)\in C(A^S)$.

Let $(i,j)\in C(A^S)$, then observe that $a_{ij}^Sy_j<y_i$ is
impossible, because the multiplication with other inequalities over
the critical cycle would lead to $1<1$. So $a_{ij}^S y_j=y_i$, and
hence there exists $k\in S_j$ such that $a_{ik} a_{jk}^{-1} y_j=y_i$.
But then also $k\in S_i$ and $(i,j)\in G_S$.
\end{proof}

The equality $G_S=C(A^S)$ means that
$C(A^S)$ is symmetrical and $(i,j)\in G_S$ if and only if $(i,j)$ or equivalently
$(j,i)$ belong to $C(A^S)$. Theorem \ref{subeig-dimension} and Proposition \ref{gs=gc}
yield the following result, see also Develin and Sturmfels \cite{DS-04}, Proposition 17.

\begin{theorem}
\label{region-dimension} Let $A\in\Rpnm$ be a matrix with no zero
rows and columns, let $y$ be a positive vector and $S$ be the row
type of $y$ with respect to $A$, then both max-algebraic and linear
dimensions of $\operatorname{cl}(X_S)$ are equal to the number of
connected components in $G_S$.
\end{theorem}

\subsection{Row and column Kleene stars}

For a matrix $A=(a_{ij})\in\Rpnn$ and any permutation $\sigma\in S_n$ (where $S_n$ denotes the group
of all permutations of $[n]$) define the {\em weight} of $\sigma$ to be
$w(\sigma):=\prod_{i=1}^n a_{i\sigma(i)}$. The {\em max-algebraic permanent} of $A$ is defined as
\begin{equation}
\label{e:perm}
\perm(A)=\bigoplus_{\sigma\in S_n} w(\sigma),
\end{equation}
and a permutation, at which the maximum in \eqref{e:perm} is attained, is called a {\em maximal
permutation}. For any permutation $\sigma$, define the diagonal matrix $D^{\sigma}=(d^{\sigma}_{ij})$ by
\begin{equation*}
d^{\sigma}_{ij}=
\begin{cases}
a_{ij}, & \text{if $j=\sigma(i)$;}\\
0, & \text{otherwise}.
\end{cases}
\end{equation*}
Observe that $A(D^{\sigma})^{-1}$ is an instance of $A^S$, for the
type $S=\{\{\sigma(1)\},\ldots,\{\sigma(n)\}\}$. The subeigencone
$V^*(A(D^{\sigma})^{-1})$ is precisely the closed region
$\operatorname{cl}(X^S)$. It contains positive vectors if and only
if $A(D^{\sigma})^{-1}$ is strongly definite, and this is true if
and only if the permutation $\sigma$ is maximal \cite{But-03}. This
is also equivalent to $(D^{\sigma})^{-1}A$ being strongly definite.
Further $A(D^{\sigma})^{-1}$ will be denoted by $A^{c\sigma}$ and
$(D^{\sigma})^{-1}A$ will be denoted by $A^{r\sigma}$. The entries
of $A^{c\sigma}$ and $A^{r\sigma}$ are
{\sloppy

}
\begin{equation}
\label{cldef}
a_{ij}^{c\sigma}=a_{i\sigma(j)}a_{j\sigma(j)}^{-1},\quad a_{ij}^{r\sigma}=
a^{-1}_{\sigma^{-1}(i) i} a_{\sigma^{-1}(i)j}.
\end{equation}
The Kleene stars of $A^{c\sigma}$ and $A^{r\sigma}$ will be denoted by $A^{c\sigma *}$ and
$A^{r\sigma *}$ and called {\em column Kleene stars} and {\em row Kleene stars}, respectively.

The results of Yoeli \cite{Yoe-61}, see also Cuninghame-Green
\cite{CG:79}, Theorem 27-11, and Izhakian \cite{Izh-05,Izh-08}
suggest that row and column Kleene stars are related to the
max-algebraic pseudoinverses of matrices. The {\em pseudoinverse} of
$A$ is defined, see \cite{CG:79} and \cite{Yoe-61}, as
$A^{\treug}=(\perm(A))^{-1} A^{\adjug}$. Here $A^{\adjug}$ is the
{\em pseudoadjugate} of $A$ defined by
$a_{ij}^{\adjug}=\perm(A_{ji})$, where $A_{ji}$ is the complementary
minor to $a_{ij}$. The following proposition collects some facts
about strongly definite matrices, which are due to Yoeli and
Cuninghame-Green.

\begin{proposition}
\label{str-def} Let $A\in\Rpnn$ be strongly definite.
\begin{itemize}
\item [1.] $I\leq A\leq A^2\leq\ldots\leq A^{n-1}=A^n=\ldots$.
\item [2.] $A^*=A^{n-1}$.
\item [3.] $A^*=A^{\adjug}=A^{\treug}$.
\end{itemize}
\end{proposition}

Izhakian \cite{Izh-05,Izh-08} studies the products $A\otimes
A^{\treug}$ and $A^{\treug}\otimes A$ over extended tropical
semiring, with the main emphasis on the questions of regularity and
rank. In this context, he proves \cite{Izh-08} that the products
$A\otimes A^{\treug}$ and $A^{\treug}\otimes A$ are Kleene stars.
Below we give an elementary proof that over max algebra, these
products are equal to column and row Kleene stars, respectively.

\begin{theorem}
\label{defclos} Let $A\in\Rpnn$ have nonzero permanent. For any
permutation $\sigma$ with maximal weight we have that $A^{c\sigma
*}=D^{\sigma} A^{\treug}=A\otimes A^{\treug}$ and $A^{r\sigma
*}=A^{\treug} D^{\sigma}=A^{\treug}\otimes A$.
\end{theorem}
\begin{proof}
Using \eqref{cldef} and the definition of $A^{\adjug}$, we write:
\begin{equation*}
\begin{split}
a_{ij}^{\adjug} &=\bigoplus_{\pi: \pi(j)=i}\prod_{k\ne j}
a_{k\pi(k)}= \bigoplus_{\pi: \pi(j)=i}\prod_{k\ne j}
a_{\sigma^{-1}\pi(k),\pi(k)} a_{k,\sigma^{-1}\pi(k)}^{c\sigma}=\\
&=\prod_{k\ne i}
a_{\sigma^{-1}(k)k}\cdot\bigoplus_{\pi:\pi(j)=i}\prod_{k\ne
j}a_{k,\sigma^{-1}\pi(k)}^{c\sigma} =\\
&=\per(A)\cdot a_{\sigma^{-1}(i)
i}^{-1}\cdot\bigoplus_{\pi:\pi(j)=\sigma^{-1}(i)}\prod_{k\ne j}
a_{k\pi(k)}^{c\sigma}=\per(A)\cdot a_{\sigma^{-1}(i) i}^{-1}
(a^{c\sigma})^{\adjug}_{\sigma^{-1}(i)j}.
\end{split}
\end{equation*}
By Proposition \ref{str-def}, $(A^{c\sigma})^{\adjug}=A^{c\sigma
*}$, so we have obtained that
$A^{\adjug}=\per(A)(D^{\sigma})^{-1}A^{c\sigma *}$, and hence
$A^{\treug}=(D^{\sigma})^{-1} (A^{c\sigma})^*$ and $D^{\sigma}
A^{\treug}=(A^{c\sigma})^*$. We now infer that
{\sloppy

}
\begin{equation*}
\begin{split}
(A\otimes A^{\treug})_{ij} &=\bigoplus_{k} a_{ik} a_{kj}^{\treug}=
\bigoplus_k a_{ik} a_{\sigma^{-1}(k) k}^{-1} a_{\sigma^{-1}(k)
j}^{c\sigma *}= \bigoplus_k a_{i\sigma^{-1}(k)}^{c\sigma}
a_{\sigma^{-1}(k)j}^{c\sigma *}=a_{ij}^{c\sigma *}.
\end{split}
\end{equation*}
Thus $A\otimes A^{\treug}=A^{c\sigma *}$. On the other hand, one can
similarly obtain that\\ $a_{ij}^{\adjug}=\per(A)
a_{i\sigma(j)}^{r\sigma *} a_{j\sigma(j)}^{-1}$ and that
$A^{\treug}\otimes A=A^{\treug} D^{\sigma}=A^{r\sigma *}$.
\end{proof}

Clearly this theorem yields the following corollary the first part of which was obtained in \cite{Ser-07}.
This corollary means that for any matrix with nonzero permanent, both
row Kleene star and column Kleene star are uniquely defined.

\begin{proposition}
\label{ser-col} Let $A\in\Rpnn$ have nonzero permanent. Then for all
permutations $\sigma$ with maximal weight, the corresponding column
Kleene stars $A^{c\sigma *}$ are equal to each other, and the row
Kleene stars $A^{r\sigma *}$ are also equal to each other.
\end{proposition}

The idea of the proof in \cite{Ser-07} was to notice that the
(sub)eigencones of $A^{c\sigma}$ are the same for all maximal
permutations $\sigma$, and to use Proposition \ref{clounique} that
any Kleene star is uniquely defined by its column span.

For a square matrix $A$, the span of its column Kleene star is the
only region of $\spann(A)$ which may have full linear dimension, and
the linear dimension of that region determines the {\em tropical
rank} of $A$, introduced by Develin et al. \cite{DSS-05}, and also
investigated by Izhakian \cite{Izh-06}. When the tropical rank is
full, the interior of span of the column Kleene star is the {\em
simple image set} of $A$ studied by Butkovi\v{c} \cite{But-00}: It
is the set of vectors $y\in\Rpn$ such that $Ax=y$ has a unique
solution. In what follows, the span of column Kleene star of $A$
will be called the {\em essential span} of $A$.

The following theorem, which is a slight generalization of Theorem 8 by Gaubert and Meunier \cite{GM-08},
illustrates the role of essential span in the geometry of max cones. It can be thought of as
a colourful generalization of Minkowski's theorem for max cones in the sense of B\'{a}r\'{a}ny
\cite{Bar-82}.

\begin{theorem}
\label{colour-mink} Let $U\subseteq \Rpn$ be a closed max cone and
let $V^1,\ldots,V^n\subseteq\Rpn$ be closed max cones such that the
intersection of $V^i$ with $U$ is nontrivial for all $i\in[n]$. Then
there exist vectors $v^1,\ldots, v^n$ such that $v^i$ is an extremal
of $V^i$, for $i\in[n]$, and $\spann(v^1,\ldots,v^n)$ has nontrivial
intersection with $U$.
\end{theorem}
\begin{proof}
Take any nonzero points $y^1\in V^1\cap U,\ldots, y^n\in V^n\cap U$
and consider the matrix $A\in\Rpnn$ with columns $A_{\cdot i}=y^i$,
for $i=1,\ldots,n$. Assume first that $A$ has permutations with
nonzero weight. The essential span of $A$ is the closed region
$\operatorname{cl}(X_S)$, where
$S=\{\{\sigma(1)\},\ldots,\{\sigma(n\}\}$, for any maximal
permutation $\sigma$. Take any $u\in\operatorname{cl}(X_S)$, then
$u\in U$ and $u\geq_i A_{\cdot \sigma(i)}$ for all $i$. The column
$A_{\cdot\sigma(i)}$ is equal to $y^{\sigma(i)}$ and it belongs to
$V^{\sigma(i)}$. Applying Minkowski theorem (Proposition \ref{mink})
and the multiorder principle (Proposition \ref{p:multiorder}), we
obtain an extremal $v^{\sigma(i)}$ of $V^{\sigma(i)}$ such that
$v^{\sigma(i)}\leq_i y^{\sigma(i)}\leq_i u$. Applying Proposition
\ref{p:multiorder} again, we see that
$u\in\spann(v^{\sigma(1)},\ldots,v^{\sigma(n)})$. As $u\in U$, the
claim follows.

In the case when $A$ does not have nonzero permutations, an inductive argument
using Hall's marriage theorem, see \cite{GM-08}, shows that there exist subsets of
indices $M$, $N_1$ and $N_2$ such that the submatrix $A_{[N_1,M]}$ is zero, while
the submatrix $A_{[N_2,M]}$ is square and has a permutation with nonzero weight.
Then the above argument goes with the essential span of that submatrix.
\end{proof}

\section{Cyclic projectors and the alternating method}

\subsection{Cyclic projectors and separation of several max cones}
Let $V^1,\ldots,V^k$ be closed max cones in $\Rpn$ and denote by
$P_i$ the projector onto $V^i$. The composition $P_k\cdots P_1$ will
be called the {\em cyclic projector} associated with $V^1,\ldots,V^k$.
This operator inherits many properties of the sole projector:
it is a homogeneous, continuous, isotone and nonincreasing operator.
In general, it is not linear with respect to max and min operations.
Such operators can be treated by nonlinear Perron-Frobenius theory.
In particular, the following theorem of Nussbaum \cite{Nus-86}
generalizes the well-known Collatz-Wielandt formula for the spectral
radius of a nonnegative matrix.
{\sloppy

}
\begin{theorem}
\label{Nussbaum} Let $F$ be a continuous, homogeneous and isotone
operator in $\Rpn$. Then the spectral radius of $F$ is equal to
\begin{equation}
\label{e:Nussbaum}
r(F)=\inf\{\lambda\mid \exists y\ \text{\rm positive:}\ Fy\leq\lambda y\}.
\end{equation}
\end{theorem}

Such operators have no more than one eigenvalue over any set of
vectors with the same support, and therefore the total number of
their eigenvalues is finite. Formula \eqref{e:Nussbaum} implies that
the spectral radius is monotone. Define the {\em cyclic projective
distance} of $y^1,\ldots,y^k\in\Rpn$ by
\begin{equation}
\label{rhdef}
\rho_{\hilb}(y^1,\ldots,y^k)=\log \bigoplus_{i_1,\ldots,i_k\in M}
y^1_{i_1} (y^2_{i_1})^{-1}\cdot\ldots\cdot y^k_{i_k}(y^1_{i_k})^{-1},
\end{equation}
when $\supp(y^1)=\ldots=\supp(y^k)=M$, and by $+\infty$ otherwise.
In the case $k=2$ this is the Hilbert projective distance between two points in $\Rpn$.
An equivalent definition is
\begin{equation}
\label{rhdef1}
\rho_{\hilb}(y^1,\ldots,y^k)=\log\inf\{\prod_{i=1}^k\lambda_i\mid y^i\leq\lambda_i y^{i+1}, i\in[k]\},
\end{equation}
where $y^{k+1}:=y^1$. Note that $\rho_{\hilb}$ is stable
under multiplication of the arguments by nonzero scalars
and under their cyclic permutation.
If $\sum_{l=1}^n y_l^i=1$ for $i\in[k]$, then it follows from \eqref{rhdef1} that
$\lambda_i\geq 1$, and $\rho_{\hilb}(y^1,\ldots,y^k)=0$
if and only if $y^1=\ldots=y^k$. For general $y^1,\ldots,y^k\in\Rpn\backslash\{0\}$,
$\rho_{\hilb}(y^1,\ldots,y^k)=0$ if and only if $y^1,\ldots,y^k$ are proportional to each other.

Define the cyclic projective distance between closed max cones
$V^1,\ldots,V^k$ by
\begin{equation}
\label{rhcdef}
\rho_{\hilb}(V^1,\ldots, V^k)=\inf_{y^1\in V^1,\ldots, y^k\in V^k} \rho_{\hilb}(y^1,\ldots y^k).
\end{equation}
The minimum in \eqref{rhcdef} is attained since $\rho_{\hilb}$ is
lower semicontinuous, see Proposition~\ref{rhos-semicont} below.

The monotonicity of spectral radius is crucial for the following
theorem \cite{GS-07}.
\begin{theorem}
\label{specrad=dh} Let $V^1,\ldots, V^k$ be closed max cones in
$\Rpn$. Suppose that $y^0$ is an eigenvector of $P_k\cdots P_1$
associated with the spectral radius, and consider vectors $y^1\in
V^1,\ldots,y^k\in V^k$ defined by $y^1:=P_1 y^0,\ldots, y^k:=P_k
y^{k-1}$. Then
\begin{equation*}
\label{specrad-cp}
\rho_{\hilb}(y^1,\ldots,y^k)=\rho_{\hilb}(V^1,\ldots,V^k)=-\log
r(P_k\cdots P_1).
\end{equation*}
\end{theorem}
Cyclic projectors also enable to prove a separation theorem for
closed max cones \cite{GS-07}, with the following ideas in mind.
Firstly, formula \eqref{e:Nussbaum} implies the existence of a
positive subeigenvector with $\lambda<1$. Secondly, if we take such
a subeigenvector, then its projections onto $V^1,\ldots,V^k$ define
separating halfspaces, see Theorem~\ref{separ0}.

\begin{theorem}
\label{separ} Let $V^1,\ldots, V^k\subseteq\Rpn$ be closed max
cones. If each of $V^1,\ldots,V^k$ has a positive vector, then the following are
equivalent.
\begin{itemize}
\item[1.] There exists a positive vector $y$ and $\lambda<1$: $P_k\cdots P_1 y\leq\lambda y$.
\item[2.] There exist halfspaces $\halfs^1,\ldots,\halfs^k$ such that
$V^1\subseteq\halfs^1,\ldots, V^k\subseteq\halfs^k$ and $\bigcap_{i=1}^k \halfs^i=\{0\}$.
\item[3.] $\bigcap_{i=1}^k V^i=\{0\}$.
\item[4.] $r(P_k\cdots P_1)<1$.
\end{itemize}
The statements 2. and 3. are equivalent even if $V^1,\ldots,V^k$ do not have
positive vectors.
\end{theorem}

\subsection{The alternating method and its convergence}

In what follows we consider the case when
$V^1=\spann(A^{(1)}),\ldots,V^k=\spann(A^{(k)})$, and
$A^{(1)},\ldots,A^{(k)}$ are nonnegative matrices with an equal number
of nonzero rows. A natural question is to find a positive solution
to the system of equations
\begin{equation}
\label{mainsystem}
A^{(1)}\otimes x^1=\ldots=A^{(k)}\otimes x^k,
\end{equation}
and the cyclic projectors provide an efficient method for doing this.

\vskip 1cm \noindent {\bf ALTERNATING METHOD}

\msn {\bf Input:} Nonnegative matrices $A^{(1)}\in\Rpnml,\ldots,
A^{(k)}\in\Rpnmk$ with an equal number $n$ of nonzero rows.

\msn {\bf Initialization:} Arbitrary positive $y^{(0)}:=y^{(1)0}$.

\msn {\bf Iteration:} Number $l\geq 1$.  For all $s=1,\ldots,k$
compute $x^{(l)s}:=\overline{A^{(s)}}\otimes' y^{(l)s-1}$ and
$y^{(l)s}:=A^{(s)}\otimes x^{(l)s}$. Set $x^{(l)}:=x^{(l)k}$ and
$y^{(l)}:=y^{(l)k}$.

\msn {\bf Stop:} If $y^{(l)}=y^{(l-1)}$, then stop. The vectors $x^{(l)s}$, for $s=1,\ldots,k$,
give a solution to system \eqref{mainsystem}. Else if $y_i^{(l)}<y_i^{(0)}$ for
all $i\in[n]$, then stop. There is no solution.

\vskip 1cm
Over the semiring
$\R_{\max,+}=(\R\cup\{-\infty\},\oplus=\max,\otimes=+)$ and for
$k=2$, this method was formulated by Cuninghame-Green and
Butkovi\v{c} \cite{CGB}. The method is essentially a max-algebraic
version of the {\em cyclic projections method} known in
optimization theory \cite{BBL-97}, since $y^{(l)}=P_k\cdots P_1
y^{(l-1)}$.

The first part of the stop condition follows from the fact that $P_1,\ldots,P_k$
are nonincreasing projectors onto $\spann(A^{(1)}),\ldots,\spann(A^{(k)})$. Indeed, if
$y^{(l-1)}=y^{(l)}$, then the
inequalities
\begin{equation*}
\label{projineqs1}
y^{(l)}\geq P_{k-1}\cdots P_1 y^{(l-1)}\geq\ldots\geq P_1 y^{(l-1)}\geq y^{(l-1)}
\end{equation*}
are satisfied with equalities, implying that $y^{(l)s}=P_s\cdots P_1
y^{(l)}$ are equal for all $s\in[k]$ and that
$y^{(l)}\in\spann(A^{(1)})\cap\ldots\cap\spann(A^{(k)})$. As
$y^{(l)s}=A^{(s)}\otimes x^{(l)s}$ for $s\in[k]$, we have that
$x^{(l)s}$, for $s\in[k]$, give a solution to \eqref{mainsystem}.

Also note that the absence of zero rows in the matrices
implies that all vectors in the sequence generated by the alternating method are positive
and hence any solution, which the alternating method may find, has to be positive.

The following proposition, similar to the results of \cite{CGB}, justifies the second
part of the stop condition. It emphasizes the role of {\em sleepers}, i.e.,
such indices
$i(s)\in [n]$ (for $s=1,\ldots,k$) that
$y^{(1)s}_{i_s)}=y^{(2)s}_{i_s}=\ldots$ for the whole sequence $\{y^{(l)s},\ l\geq 1\}$,
and $j(s)\in [m_s]$ such that $x^{(1)s}_{j_s}=x^{(2)s}_{j_s}=\ldots$
for the whole sequence $\{x^{(l)s},\ l\geq 1\}$.
Sleepers will be called {\em eternal}, if the corresponding coordinates
are constant for all $l\geq 1$, and {\em temporary}, if the corresponding coordinates
are constant up to the last iteration of the alternating method.

\begin{proposition}
\label{sleeper}
Let $A^{(1)},\ldots,A^{(k)}$, $x^{(l)s}$ and $y^{(l)s}$ be as in the
formulation of the alternating method. Then
\begin{itemize}
\item[1.] temporary sleepers exist for
all sequences $\{x^{(l)s}\}$ and $\{y^{(l)s}\}$, $s\in [k]$.
\item[2.] if \eqref{mainsystem} has a solution, then
eternal sleepers exist for all sequences
$\{x^{(l)s}\}$ and $\{y^{(l)s}\}$, $s\in [k]$.
\item[3.] if \eqref{mainsystem} has a positive solution, then
$\{x^{(l)s}\}$ and $\{y^{(l)s}\}$, for all $s\in[k]$, are bounded from below by positive vectors.
\end{itemize}
\end{proposition}
\begin{proof}

1. Assume that for some $s\in[k]$ and $l\geq 1$ we have that all
coordinates of $y^{(l)s}$ or $x^{(l)s}$ are strictly less than that
of $y^{(1)s}$ or $x^{(1)s}$. Then we have that $y^{(l)s}\leq\mu
y^{(1)s}$ or $x^{(l)s}\leq\mu x^{(1)s}$ for some $\mu<1$. As all
matrix multiplications are homogeneous and isotone, we have that
$y^{(l)}\leq\mu y^{(1)}$ so that all coordinates of $y^{(l)}$ are
strictly less than that of $y^{(0)}$ and the alternating method
immediately stops.

2. and 3. Take any $s\in[k]$. If there is a vector $y$ in the
intersection of column spans, we can scale it so that $y\leq
y^{(1)s}$ and $y_i=y_i^{(1)s}$ for some $i$. In terms of the
multiorder, $y\leq_i y^{(1)s}$ (for this scaling it is essential
that $y^{(0)}$ and hence $y^{(1)s}$ are positive). As the projectors
are all isotone and $y$ is their fixed point, we have that $y\leq
y^{(l)s}$ and $y_i=y_i^{(l)s}$ for the whole sequence. If
\eqref{mainsystem} has a positive solution, then the same scaling
argument shows that the sequence $\{y^{(1)s},y^{(2)s},\ldots\}$ is
bounded from below by a positive vector. Now note that the same line
of argument applies to $\{x^{(l)s}\}$ as well.
\end{proof}

In what follows we will prove that the alternating method converges to a positive solution if a positive
solution exists. We note here that a cyclic projector
is a min-max function in the sense of \cite{CZ-02,CTGG-99,OP-97}, 
with addition being replaced by multiplication, 
and the convergence of the alternating method
follows from the results of \cite{CZ-02,OP-97} concerning the ultimate periodicity
of min-max functions. Below we give a different proof which uses the cellular decomposition idea.

We first investigate the convergence of the alternating method for Kleene stars,
which then enables us, using cellular decomposition, to prove the finiteness
results for general matrices.

\begin{proposition}
\label{poly-k} Suppose that $A^{(1)},\ldots,A^{(k)}\in\Rpnn$
are Kleene stars. If $\spann(A^{(1)})\cap\ldots \cap\spann(A^{(k)})$
contains a positive vector, then the alternating method converges in
no more than $n$ iterations.
\end{proposition}
\begin{proof}
The alternating method starts with an arbitrary positive initial
vector $y$ and repeatedly applies the composition $P_k\cdots P_1$.
Due to Proposition \ref{proj-linear} we have that
\begin{equation*}
P_k\cdots P_1 y=\overline{A^{(k)}}\otimes'\ldots\otimes'\overline{A^{(1)}}\otimes' y,
\end{equation*}
and hence
\begin{equation*}
(P_k\cdots P_1)^m y=(\overline{A^{(k)}}\otimes'\ldots\otimes'\overline{A^{(1)}})^m\otimes' y.
\end{equation*}
This means that the stabilization of the alternating method is
equivalent to the stabilization of $(A^{(1)}\otimes\ldots\otimes
A^{(k)})^m\otimes y$ for any positive $y$. Denote the matrix product
$A^{(1)}\otimes\ldots\otimes A^{(k)}$ by $C$. By Proposition
\ref{star-inters-k} we have that $\lambda(C)=1$. We also have that
the diagonal entries of $C$ are equal to $1$ and hence it is a
strongly definite matrix. By Proposition \ref{str-def} the powers of
$C$ stabilize in no more than $n-1$ steps, and this proves the
claim.
\end{proof}

Now we make use of the cellular decomposition to prove
that if there is a positive solution, then the alternating method finds a positive
solution in a finite number of steps. First we prove the following technical proposition.

\begin{proposition}
\label{liberal} Suppose that $A^{(1)},\ldots,A^{(k)}\in\Rpnn$ have
all diagonal entries equal to $1$ and suppose that any product $D$
of no more than $n$ of them has $\lambda(D)\leq 1$. Fix a mapping
$j:\{1,\ldots\}\mapsto\{1,\ldots,k\}$. Consider the sequence of
products $C^{(m)}=A^{(j(m))}\otimes\ldots\otimes A^{(j(1))}$, for
$m\geq 1$. Then there exists $m\leq n^k-1$ such that
$C^{(m)}=C^{(m+1)}$.
\end{proposition}

\begin{proof}
For the case of just one matrix, this is Proposition \ref{str-def}.
We argue by induction, assuming the result is true for $k-1$
matrices and proving it for $k$. Choose any mapping $\pi
:\{1,\ldots,n\}\mapsto\{1,\ldots,k\}.$ Then either for some $m<n^k$
we have that there are no repetitions before that $m$ and $$
C^{(m)}=\bigotimes_{i=1}^n A^{(\pi(i))}\otimes  B^{(i)},$$ where
each $B^{(i)}$ is a product of less than $n^{k-1}-1$ matrices, or
there is a repetition, and in this case we are done. Hence, for
$M=n^k-1$, either there are repetitions before that $M$, or the
product $C^{(M)}=(c^{(M)}_{ij})$ contains all the above mentioned
products. We claim then that
\begin{equation}
\label{cij-rep} c^{(m)}_{ij}=\bigoplus_{\pi, i_{n-1},\ldots, i_1}
a_{i\, i_{n-1}}^{(\pi(n))}\cdot\ldots\cdot a_{i_1\, j}^{(\pi(1))}.
\end{equation}
for all $m\geq M$. Indeed, $c^{(m)}_{ij}$ is greater than or equal
to the maximum on the r.h.s. due to the choice of $M$ and since all
diagonal entries of all matrices are $1$. It is actually equal to
this maximum because all products of no more than $n$ matrices have
$\lambda\leq 1$, so the weight of any path of length $M$ does not
exceed the weight of the simple path obtained after cycle deletion,
and the weights of all simple paths are already in \eqref{cij-rep}.
\end{proof}

\begin{theorem}
\label{finite}
Suppose that $A^{(1)}\in\Rpnml,\ldots,A^{(k)}\in\Rpnmk$ have all rows nonzero and are such that
$\spann(A^{(1)})\cap\ldots\cap\spann(A^{(k)})$
contains a positive vector. Then the alternating method stabilizes in a finite number of steps.
\end{theorem}

\begin{proof}
It follows from Proposition \ref{cellular2} that  for each matrix
$A^{(i)}$ we have a Kleene decomposition
\begin{equation*}
\spann(A^{(i)})=\bigcup_{l=1}^{s(i)} \spann(A^{(il)}),
\end{equation*}
where $A^{(il)}\in\Rpnn$ are Kleene stars. Then we have that

\begin{equation}
\label{papbm}
\begin{split}
&(P_k\cdots P_1)^m y=(\overline{A^{(kl(k,m))}}\otimes'\ldots\otimes'\overline{A^{(1l(1,m))}})
\otimes'\ldots\\
&\otimes'(\overline{A^{(kl(k,1))}}\otimes'\ldots\otimes'\overline{A^{(1l(1,1))}})\otimes' y
\end{split}
\end{equation}
for some index mappings $l(i,j).$

\msn It suffices to prove the
stabilization of the sequence
\begin{equation}
\label{aimbjm}
\begin{split}
B^{(m)}\otimes\ldots\otimes B^{(1)} \otimes y,
\end{split}
\end{equation}
where $B^{(i)}=((A^{(kl(k,i))})^T\otimes\ldots\otimes
(A^{(1l(1,i))})^T)$. Note that the number of matrices $B^{(i)}$ is
also finite. Since the spans of the matrices
$A^{(1)},\ldots,A^{(k)}$ have a point in intersection, by
Proposition \ref{sleeper} sequence \eqref{papbm} is bounded from
below, and hence \eqref{aimbjm} is bounded from above.

\msn Consider a finite product $B$ of some matrices $B^{(i)}$, appearing in \eqref{aimbjm}.
%and taken in arbitrary order.
If $\lambda(B)>1$, then at least one
of the matrices making this product
will appear only a finite number of times. Otherwise the sequence will be unbounded,
which is a contradiction.

\msn Hence after some finite $m$ the matrices $B^{(i)}$ appearing in
the sequence will be such that $\lambda(B)\leq 1$ for any product
$B$ of no more than $n$ of them.

\msn After that, the finite convergence of alternating method is
guaranteed by Proposition \ref{liberal}.
\end{proof}

\subsection{Bounds on the number of iterations}

Now we examine the case when the system has no solution, i.e., when
the max cones $\spann(A^{(1)}),\ldots,\spann(A^{(k)})$ do not have nontrivial
intersection. Here we will need the
{\em total projective distance} between $y^1,\ldots,y^k$, which is the sum of
projective distances
\begin{equation}
\label{totaldist-vecs}
\rho_{\Sigma}(y^1,\ldots,y^k)=\rho_{\hilb}(y^1,y^2)+\ldots +\rho_{\hilb}(y^k,y^1),
\end{equation}
if $y^1,\ldots,y^k$ have equal supports, and $+\infty$ otherwise. Note that
\begin{equation}
\label{cycl-total}
\rho_{\Sigma}(y^1,\ldots,y^k)=\rho_{\hilb}(y^1,\ldots,y^k)+\rho_{\hilb}(y^k,\ldots,y^1),
\end{equation}
where $\rho_{\Hilb}$ is the cyclic projective distance defined by \eqref{rhdef}.
By analogy with \eqref{rhdef1},
\begin{equation}
\label{totaldist1}
\rho_{\Sigma}(y^1,\ldots,y^k)=\log\inf\{\prod_{i=1}^k \lambda_i\mu_i\mid y^i\leq\lambda_i y^{i+1},\,
y^{i+1}\leq\mu_i y^i,\, i\in[k]\},
\end{equation}
where $y^{k+1}:=y^1$. Like $\rho_{\hilb}$, the total projective distance is stable under
scalar multiplication of the arguments and their cyclic permutation.

Denote $\cS_n:=\{x\in\Rpn\mid \sum_{i=1}^n x_i=1\}$ and consider
$\cS_n^k:=\overbrace{\cS_n\times\ldots\times\cS_n}^{k}$ endowed with product topology.
A function $\phi: \cS_n^k\mapsto\Rp\cup\{+\infty\}$ is called {\em lower semicontinuous} if the sublevel sets
\begin{equation}
\cS_n^k(\phi,a)=\{(y^1,\ldots,y^k)\in\cS_n^k\mid\phi(y^1,\ldots,y^k)\leq a\},
\end{equation}
are closed for all $a\in\Rp$. The author gratefully acknowledges the idea of the proof of the following proposition to
St\'{e}phane Gaubert.
\begin{proposition}
\label{rhos-semicont}
$\rho_{\Sigma}(y^1,\ldots,y^k)$ and $\rho_{\hilb}(y^1,\ldots,y^k)$
are lower semicontinuous on $\cS_n^k$.
\end{proposition}
\begin{proof}
Consider sequences $\{y^{(m)i},\ m\geq 1\}\subseteq\cS_n$ converging
to $\overline{y^i}$, for $i\in[k]$. We need to show that if
$(y^{(m)1},\ldots,y^{(m)k})\in\cS_n^k(\rho_{\Sigma},a)$ (resp. if
$(y^{(m)1},\ldots,y^{(m)k})\in\cS_n^k(\rho_{\hilb},a)$) for all $m\geq 1$,
then $(\overline{y^{1}},\ldots,\overline{y^{k}})\in\cS_n^k(\rho_{\Sigma},a)$
(resp. $(\overline{y^{1}},\ldots,\overline{y^{k}})\in\cS_n^k(\rho_{\hilb},a)$).
If $(y^{(m)1},\ldots,y^{(m)k})\in\cS_n^k(\rho_{\Sigma},a)$ for all $m$, there exist
$\lambda_i^{(m)},\mu_i^{(m)}\in\Rp$ such that $y^{(m)i}\leq\lambda_i^{(m)}y^{(m)i+1}$
and $y^{(m)i+1}\leq\mu_i^{(m)}y^{(m)i}$ for all $i\in[k]$,
and that $\prod_{i=1}^k\lambda_i^{(m)}\mu_i^{(m)}\leq a$.
As $\sum_{l=1}^n y^{(m)i}_l=1$ for all $m$
and $i$, and $y^{(m)i}\leq\lambda_i^{(m)}y^{(m)i+1}$
and $y^{(m)i+1}\leq\mu_i^{(m)}y^{(m)i}$, we have that $\lambda_i^{(m)}\geq 1$ and $\mu_i^{(m)}\geq 1$.
Using these inequalities and $\prod_{i=1}^k\lambda_i^{(m)}\mu_i^{(m)}\leq a$,
we obtain that $1\leq\lambda_i^{(m)}\leq a$ and $1\leq\mu_i^{(m)}\leq a$ for all $i\in[k]$.
Taking convergent subsequences if necessary, we can assume that $\lambda_i^{(m)}\to\lambda_i$
and $\mu_i^{(m)}\to\mu_i$
for $i\in[k]$. Then we have $\overline{y^i}\leq\lambda_i\overline{y^{i+1}}$
and $\overline{y^{i+1}}\leq\mu_i\overline{y^i}$
for all $i\in[k]$, and $\prod_{i=1}^k \lambda_i\mu_i\leq a$, which yields
$(\overline{y^1},\ldots,\overline{y^k})\in\cS_n^k(\rho_{\Sigma},a)$. The proof for the case
of $\rho_{\Sigma}$ is complete, the case of $\rho_{\hilb}$ is treated analogously.
\end{proof}

By analogy with \eqref{rhcdef},
the total projective distance between closed max cones $V^1,\ldots,V^k$ is defined by
\begin{equation}
\label{totaldist-cones}
\begin{split}
\rho_{\Sigma}(V^1,\ldots,V^k)&=\rho_{\hilb}(V^1,V^2)+\ldots+\rho_{\hilb}(V^k,V^1)=\\
&=\min_{y^1\in V^1,\ldots,y^k\in V^k} \rho_{\Sigma}(y^1,\ldots, y^k).
\end{split}
\end{equation}
Observe that $\rho_{\Sigma}(y^1,\ldots,y^k)=0$
if and only if $y^1,\ldots,y^k$ are multiples of each other. This is generalised in the
following proposition.
\begin{proposition}
\label{p:sigmahilb}
Let $V^1,\ldots, V^k\subseteq\Rpn$ be closed max cones. Then
$\rho_{\Sigma}(V^1,\ldots, V^k)=0$ (equivalently, $\rho_{\hilb}(V^1,\ldots,V^k)=0$)
if and only if the intersection
of $V^1,\ldots, V^k$ is nontrivial.
{\sloppy

}
\end{proposition}
\begin{proof}
We show the ``only if'' part. The intersections of $V^i$ and $\cS_n$
are closed sets. Let the sequences $\{y^{(m)i},\ m\geq 1\}$, for
$i\in[k]$ and $y^{(m)i}\in V^i\cap \cS_n$, be such that
$\lim_{m\to\infty}\rho_{\Sigma}(y^{(m)1},\ldots, y^{(m)k})=0$ (or
$\lim_{m\to\infty}\rho_{\hilb}(y^{(m)1},\ldots, y^{(m)k})=0$). As
$\cS_n$ is compact, we can assume that $y^{(m)i}\to\overline{y^i}$
for $i\in[k]$, where $\overline{y^i}\in V^i\cap\cS_n$ as
$V^i\cap\cS_n$ is closed. Proposition \ref{rhos-semicont} implies
that $\rho_{\Sigma}(\overline{y^1},\ldots,\overline{y^k})=0$ (resp.
$\rho_{\hilb}(\overline{y^1},\ldots,\overline{y^k})=0$). Hence
$\overline{y^i},$ for $i\in[k]$, are proportional vectors contained
in $V^1\cap\ldots\cap V^k$. The proof of the ``only if'' part is
complete. The ``if'' part is obvious.
\end{proof}

Let vector $y$ and matrix $A$ have finite entries. Denote
\begin{equation}
\label{lognorms}
||y||=\log\bigoplus_{i,j} y_iy_j^{-1},\quad ||A||=\log\bigoplus_{i,j,k} a_{ik} a_{jk}^{-1}.
\end{equation}
A vector $y=\bigwedge_{i=1}^n \lambda_i A_{\cdot i},$ where $\lambda_i>0$ for all $i\in [n]$,
and $\wedge$ denotes the componentwise minimum, will be called a {\em min combination}
of the columns of $A$.
\begin{proposition}
\label{norms}
Let $A\in\Rpnm$ and $y\in\Rpn$ have all entries positive. If $y$ is a max combination
or a min combination of the columns of $A$, then $||y||\leq||A||$.
\end{proposition}
\begin{proof}
Let $y=\bigoplus_j \lambda_j A_{\cdot j}$, or let $y=\bigwedge_j \lambda_j A_{\cdot j}$ with
all $\lambda_j\neq 0$. Then
\begin{equation*}
\label{norms-est}
\begin{split}
\exp(||y||)&=\bigoplus_{i,j} y_iy_j^{-1}=\bigoplus_{i,j}(\bigoplus_k \lambda_k a_{ik})\cdot
(\bigwedge_{l:\lambda_l\neq 0}\lambda_l^{-1} a_{jl}^{-1})=\\
&=\bigoplus_{i,j,k}\lambda_k a_{ik}\cdot
(\bigwedge_{l:\lambda_l\neq 0}\lambda_l^{-1} a_{jl}^{-1})\leq\bigoplus_{i,j,k:\lambda_k\neq 0}
a_{ik} a_{jk}^{-1}\leq \exp(||A||),\ \text{or}\\
\exp(||y||)&=\bigoplus_{i,j} y_iy_j^{-1}=\bigoplus_{i,j}(\bigwedge_k \lambda_k a_{ik})\cdot
(\bigoplus_l\lambda_l^{-1} a_{jl}^{-1})=\\
&=\bigoplus_{i,j,l}\lambda_l^{-1} a_{jl}^{-1}\cdot
(\bigwedge_k\lambda_k a_{ik})\leq\bigoplus_{i,j,l}
a_{il} a_{jl}^{-1}\leq\exp(||A||),
\end{split}
\end{equation*}
respectively. The claim follows by the monotonicity of the logarithm.
\end{proof}

\begin{proposition}
\label{simplest}
Let $u\in\Rpn$ be a positive vector, let $V\subseteq\Rpn$ be a closed max cone
and let $v=P_V(u)$. Then $\sum_{i=1}^n (\log u_i-\log v_i)\geq\rho_{\hilb}(u,v)$.
\end{proposition}
\begin{proof}
As $v\leq u$ and $u_k=v_k$ for some $k$, we have that $\rho_{\hilb}(u,v)=\max_{i=1}^n(\log u_i-\log v_i)$.
As any sum of nonnegative numbers is greater than or equal to any of its
terms, the claim follows.
{\sloppy

}
\end{proof}

\begin{proposition}
\label{forx} Suppose that $A\in\Rpnm$, and suppose that
$x^1,x^2\in\Rpm$ and $y^1,y^2\in\Rpn$ are positive and such that
$y^1\geq y^2$ with strict inequalities in at most $n'$
coordinates, $x^1\geq x^2$ and $A\otimes x^1=y^1$, $A\otimes
x^2=y^2$. Then
\begin{itemize}
\item[1.] there exists $k$ such that $x_k^1(x_k^2)^{-1}\geq\max_s y_s^1 (y_s^2)^{-1}$;
\item[2.] the inequality
%\begin{equation*}
%\label{total01}
$\sum_{k=1}^m (\log x_k^1-\log x_k^2)\geq \frac{1}{n'} \sum_{i=1}^n
(\log y_i^1-\log y_i^2)$ holds.
%\end{equation*}
\end{itemize}
\end{proposition}
\begin{proof}
Let $t$ be such that $\max_s y_s^1 (y_s^2)^{-1}=y_t^1 (y_t^2)^{-1}$ and define $k$
such that $\max_s (a_{ts} x_s^1)=a_{tk}x_k^1=y_t^1$.
The inequalities $a_{tk}\neq 0$ and $a_{tk}x_k^2\leq y_t^2$ imply part 1. To obtain part 2. we recall that
any sum of nonnegative numbers is greater than or equal to any of its terms, and that the maximum
is always greater than or equal to the arithmetic mean.
\end{proof}

Now we obtain a bound for the number of iterations of the alternating method.
For brevity, we denote
$\rho_{\Sigma}(A^{(1)},\ldots,A^{(k)}):=\rho_{\Sigma}(\spann(A^{(1)}),\ldots,\spann(A^{(k)}))$.

\begin{theorem}
\label{estimate1} Suppose that
$A^{(1)}\in\Rpnml,\ldots,A^{(k)}\in\Rpnmk$, that $A^{(k)}$ has all
entries positive, and that
$\spann(A^{(1)})\cap\ldots\cap\spann(A^{(k)})=\{0\}$. Then after no
more than
\begin{equation}
\label{e:estimate11}
2(n-1)\min(||A^{(k)}||,(m_k-1)||A^{(k)T}||)/\rho_{\Sigma}(A^{(1)},\ldots,A^{(k)})
\end{equation}
iterations the alternating method will terminate.
\end{theorem}

\begin{proof}
Let the sequences $\{y^{(l)s},\, l\geq 1\}$ and $\{x^{(l)s},\, l\geq
1\}$, for $s\in[k]$, be as in the formulation of the alternating
method. Using Proposition~\ref{simplest}, we obtain the following
lower bound for the total sum of logarithmic coordinate losses of
$y^{(l)}$ at each iteration:
\begin{equation}
\label{e:total1}
\begin{split}
\sum_{i=1}^n(\log y_i^{(l+1)}-\log y_i^{(l)})&=\sum_{s=0}^{k-1}\sum_{i_s=1}^n
(\log y^{(l)s+1}_{i_s}-\log y^{(l)s}_{i_s})\geq\\
&\geq\rho_{\Sigma}(y^{(l)1},\ldots,y^{(l)k})\geq\rho_{\Sigma}(A^{(1)},\ldots,A^{(k)}).
\end{split}
\end{equation}
Using Proposition~\ref{forx}, we also obtain that
\begin{equation}
\label{e:total2}
\begin{split}
\sum_{i=1}^n(\log x_i^{(l+1)}-\log x_i^{(l)})&\geq\frac{1}{n-1}\sum_{i=1}^n
(\log y_i^{(l+1)}-\log y_i^{(l)})\geq\\
&\geq\frac{1}{n-1}\rho_{\Sigma}(A^{(1)},\ldots,A^{(k)}).
\end{split}
\end{equation}
Let $j$ be a temporary sleeper for $\{x^{(l)}\}$ and let $i$ be a
temporary sleeper for $\{y^{(l)}\}$. The existence of temporary
sleepers was shown in Proposition~\ref{sleeper}. Thus the total sum
of all logarithmic coordinate losses of $y^{(l)}$ at each iteration
is at least $\rho_{\Sigma}(A^{(1)},\ldots,A^{(k)})$, while the $i$th
coordinate of $y^{(l)}$ is a sleeper, and the total sum of all
logarithmic coordinate losses of $x^{(l)}$ is at least
$\frac{1}{n-1}\rho_{\Sigma}(A^{(1)},\ldots,A^{(k)})$ while the $j$th
coordinate of $x^{(l)}$ is a sleeper. This will stop the alternating
method. Indeed, we repeatedly apply $P_k\cdots P_1$ and stop when
all coordinates of $y^{(l)}$ decrease with respect to that of
$y^{(0)}$. As $y^{(l)}$, for $l\geq 1$, is a max combination of the
columns of $A^{(k)}$, by Proposition \ref{norms} we have that $\log
y^{(1)}_t-\log y^{(1)}_i\leq||y^{(1)}||\leq||A^k||$ for all
$t\in[n]$. Lower bound \eqref{e:total1} for the total sum of
logarithmic coordinate losses of $y^{(l)}$ at each iteration implies
that after at most
$2(n-1)||A^{(k)}||/\rho_{\Sigma}(A^{(1)},\ldots,A^{(k)})$ iterations
there will be $t$ such that $\log y^{(l)}_i-\log
y^{(l)}_t>||A^{(k)}||$, if the method does not stop, and this
contradicts Proposition \ref{norms}. Hence, after at most that
number of iterations all coordinates will have to fall in value with
respect to the coordinates of the initial vector. Now, as $x^{(l)}$,
for $l\geq 1$, is a min combination of the columns of
$\overline{A^{(k)}}$, by Proposition \ref{norms} we have that $\log
x^{(1)}_t-\log x^{(1)}_i\leq||x^{(1)}||\leq||A^{(k)T}||$ for all
$t\in[m_k]$ (note that $||\overline{A}||=||A^T||$ for any positive
matrix $A$). Using \eqref{e:total2} instead of \eqref{e:total1} and
arguing as above, we obtain the upper bound
$2(m_k-1)||A^{(k)T}||/(\frac{1}{n-1}\rho_{\Sigma}(A^{(1)},\ldots,A^{(k)}))$
on the number of iterations, and this proves the claim.
\end{proof}

If there is more than one matrix with all entries positive, then bound
\eqref{e:estimate11} can be improved.

\begin{theorem}
\label{estimate2} Suppose that
$A^{(1)}\in\Rpnml,\ldots,A^{(k)}\in\Rpnmk$, that
$A^{(r_1)},\ldots,A^{(r_s)}$ have all entries positive, and that
$\spann(A^{(1)})\cap\ldots\cap\spann(A^{(k)})=\{0\}$. Then after no
more than
{\sloppy

}
\begin{equation}
\label{e:estimate22}
2(n-1)\min_{i=1}^s \min(||A^{(r_i)}||,(m_{r_i}-1)||A^{(r_i)T}||)/
\rho_{\Sigma}(A^{(1)},\ldots,A^{(k)})
\end{equation}
iterations the
alternating method will terminate.
\end{theorem}
\begin{proof}
Applying the argument of Theorem \ref{estimate1} and using the fact
that $\rho_{\Sigma}$, like $\rho_{\Hilb}$, is stable under the
cyclic permutations of its arguments, we obtain that for any
$t=1,\ldots,s$, after at most
\begin{equation}
\label{e:intermed}
l=2(n-1)\min (||A^{(r_t)}||,(m_{r_t}-1)||A^{(r_t)T}||)/\rho_{\Sigma}(A^{(1)},\ldots,A^{(k)}))
\end{equation}
iterations all coordinates of
$y^{(l)\,r_t}$ have to fall with respect to the coordinates
of $y^{(1)\,r_t}$. This means that there is a $\mu<1$ such that
$y^{(l)\,r_t}\leq\mu y^{(1)\,r_t}$. As all projectors are homogeneous and order preserving,
we also have that $y^{(l)}\leq\mu y^{(1)}$. Therefore
all the coordinates of $y^{(l)}$ decrease
with respect to that of $y^{(1)}$, and hence to that of $y^{(0)}$, and the alternating
method stops with negative answer. So the number of iterations does not exceed
\eqref{e:intermed} for each $r_t$, and hence it does not
exceed the minimum of these, which is \eqref{e:estimate22}.
\end{proof}

Now we show that the techniques developed above apply to the case of integer matrices
over the max-plus semiring $\R_{\max,+}=(\R\cup\{-\infty\},\oplus=\max,\otimes=+)$
investigated by Cuninghame-Green and Butkovi\v{c} \cite{CGB}. In what follows, we switch to
the matrix algebra over the max-plus semiring and to the alternating method formulated
over that semiring.

First note that if $y\in\Rn$ is a max-plus or min-plus combination
of columns of a matrix $A\in\Rnm$ with real entries, then $||y||\leq
||A||$, where like in \eqref{lognorms} but without logarithm, the
norms are defined by
\begin{equation}
\label{realnorms} ||y||=\max_{i,j}(y_i-y_j),\quad ||A||=\max_{i,j,k}
(a_{ik}-a_{jk}).
\end{equation}

\begin{theorem}
\label{altmeth-integer} Suppose that $A^{(1)}\in\R^{n\times
m_1},\ldots,A^{(k)}\in\R^{n\times m_k}$ have all entries integer.
Then after no more than
\begin{equation}
\label{e:estimate33}
2\min_{i=1}^k\min((n-1)\frac{k-1}{k}||A^{(i)}||,(m_i-1)||A^{(i)T}||)
\end{equation}
iterations the alternating method will terminate.
\end{theorem}
\begin{proof}
We are in almost the same situation as in Theorem \ref{estimate2}: for all
$x^{(l)s}$ and $y^{(l)s}$ there exist temporary sleepers, the norms
$||y^{(l)s}||$ do not exceed $||A^{(s)}||$ and the norms $||x^{(l)s}||$ do not exceed $||A^{(s)T}||$.
It remains to give bounds for the total sum of coordinate losses for $x^{(l)s}$ and $y^{(l)s}$ at
each iteration. As everything is integer, the total sum of losses for both
$x^{(l)s}$ and $y^{(l)s}$ is not
less than $1$. The multiple $\frac{k-1}{k}$ at $||A^{(i)}||$, which may be important
only if $k$ is small, is due to the
observation that if we apply $P_1,\ldots,P_{k-1}$ to $y^{(l)}\in A^{(k)}$ and do
not see any fall in coordinates,
then $y^{(l)}$ is in the intersection and the method immediately stops, hence
during the run of the algorithm, after at most $k-1$ actions (not $k$ but $k-1$) of the sole projectors
at least one coordinate of $y$ has to fall.
The claim now follows by the same argument as in Theorems \ref{estimate1} and \ref{estimate2}.
\end{proof}

The bounds on number of iterations in \cite{CGB}, obtained in the case $k=2$,
are in the same vein as \eqref{e:estimate33}. The only bound on number of iterations in \cite{CGB}
which does not depend on the choice of initial vector
would read in our terms essentially as
$2\min_{l=1}^k ((m_l-1)\max_{i,j}(|a^{(l)}_{ij}|))$, where $|\cdot|$
denotes the modulus of an entry. The bound of \eqref{e:estimate33} is expressed
in terms of projective norms of rows and columns of the matrices, which makes it more precise.

\section{Acknowledgement}
The author is grateful to the anynomous referee for careful reading and useful remarks,
and wishes to thank
%Peter Butkovi\v{c} for numerous discussions and suggestions which helped
%to improve this paper, St\'{e}phane Gaubert for his ideas concerning the proofs, in particular,
%of Proposition \ref{p:sigmahilb}
Peter Butkovi\v{c}, St\'{e}phane Gaubert,
Alexander Guterman, Grigori Litvinov, Victor Maslov, Alexander
Mikhalev, Hans Schneider and Andre\u{\i} Sobolevski\u{\i}, whose
help, interest, valuable ideas and stimulating discussions have been
very important for this work.

\bibliographystyle{amsplain}
\providecommand{\bysame}{\leavevmode\hbox to3em{\hrulefill}\thinspace}
\providecommand{\MR}{\relax\ifhmode\unskip\space\fi MR }
\providecommand{\MRhref}[2]{%
  \href{http://www.ams.org/mathscinet-getitem?mr=#1}{#2}
}
\providecommand{\href}[2]{#2}

\end{document}